\documentclass[smallextended,numbook,runningheads]{svjour3}     % onecolumn (second format)
\smartqed  % flush right qed marks, e.g. at end of proof
\usepackage{amssymb,amsfonts,amsmath}
\usepackage[mathscr]{eucal}
\usepackage{stmaryrd} % \boxast
\usepackage{bm}
\usepackage{algorithm,algpseudocode} % before cleveref
\usepackage{xparse}
\usepackage[caption=false]{subfig}
\usepackage{booktabs} % Three-line table
\usepackage{graphicx} % For \resizebox
\usepackage{url}
\usepackage{hyperref}
\usepackage[nameinlink]{cleveref} % \Cref{}

\DeclareMathOperator{\Trace}{Trace}
\DeclareMathOperator{\TR}{TR}
\DeclareMathOperator{\rank}{rank}

\NewDocumentCommand \tensor {O{}m} {\boldsymbol{#1\mathscr{\MakeUppercase{#2}}}} % using package {eucal} with option [mathscr]
\newcommand{\mat}[1]{\mathbf{#1}}
\newcommand{\vect}[1]{\bm{#1}}
\newcommand{\mc}[1]{\mathcal{#1}}
\newcommand{\mf}[1]{\mathfrak{#1}}
\newcommand{\bb}[1]{\mathbb{#1}}
\newcommand{\bigO}[1]{\mathcal{O}\left( #1 \right)}
\newcommand{\fn}{\mc{F}_n}
\newcommand{\hn}{\mc{H}_n}
\newcommand{\fxi}{\mc{F}_{\xi^{t}}}
\newcommand{\hxi}{\mc{H}_{\xi^{t}}}

\newcommand{\tabincell}[2]{\begin{tabular}{@{}#1@{}}#2\end{tabular}}
\catcode`\|=12
\begin{document}
\begin{sloppypar}

\title{Block-Randomized Stochastic Methods for Tensor Ring Decomposition\thanks{The work is supported by the National Natural Science Foundation of China (No. 11671060) and the
Natural Science Foundation of Chongqing, China (No. cstc2019jcyj-msxmX0267). 
% General acknowledgments should be placed at the end of the article.
}}
% \subtitle{Do you have a subtitle?\\ If so, write it here}


\author{Yajie Yu         \and
        Hanyu Li         \and    %etc. 
        Jingchun Zhou
}


\institute{Yajie Yu \at
              College of Mathematics and Statistics, Chongqing University, Chongqing, 401331, P.R. China;
              % Tel.: +123-45-678910\\
              % Fax: +123-45-678910\\
              \email{zqyu@cqu.edu.cn}           %  \\
           \and
           Hanyu Li \at
              College of Mathematics and Statistics, Chongqing University, Chongqing, 401331, P.R. China;
              \email{lihy.hy@gmail.com or hyli@cqu.edu.cn}
           \and
           Jingchun Zhou \at
              College of Mathematics and Statistics, Chongqing University, Chongqing, 401331, P.R. China;
              \email{zhoujc@cqu.edu.cn}
}

\date{Received: date / Accepted: date}

\maketitle

\begin{abstract}
Tensor ring (TR) decomposition is a simple but effective tensor network for analyzing and interpreting latent patterns of tensors.
In this work, we propose a doubly randomized optimization framework for computing TR decomposition. 
It can be regarded as a sensible mix of randomized block coordinate descent and stochastic gradient descent, and hence functions in a double-random manner and can achieve lightweight updates and a small memory footprint.
Further,  to improve the convergence, especially for ill-conditioned problems, we propose a scaled version of the framework that can be viewed as an adaptive preconditioned or diagonally-scaled variant.
Four different probability distributions for selecting the mini-batch and the adaptive strategy for determining the step size are also provided. 
Finally, we present the theoretical properties and numerical performance for our proposals.
%Include keywords and mathematical subject classification numbers as needed.
\keywords{Tensor ring decomposition \and Random sampling \and Stochastic gradient descent \and Scaled gradient descent \and Adaptive algorithm \and Randomized algorithm}
\subclass{15A69 \and 68W20 \and 90C52}
\end{abstract}

%%=================%%
%%  Introduction   %%
%%=================%%
\section{Introduction}
\label{sec:introduction}
Tensor, which provides a potent model for representing multi-way data, plays an indispensable role in contemporary data science with ubiquitous applications, notably for modeling high-dimensional functions or operators. 
Tensor decompositions, due to their high compression and data representation abilities, have been investigated extensively and applied widely in a variety of fields such as image and video completion, signal processing, machine learning, chemometrics, and so on; 
see the detailed review in \cite{kolda2009TensorDecompositions,sidiropoulos2017TensorDecomposition}. 
CANDECOMP/PARAFAC (CP) decomposition and Tucker decomposition are the two most well-known and traditional tensor decomposition models \cite{kolda2009TensorDecompositions}. 
Nevertheless, due to the NP-hard of finding CP decomposition and the curse of dimensionality of Tucker decomposition, tensor network decomposition has been proposed, with tensor train (TT) decomposition \cite{oseledets2011TensorTrainDecomposition} and tensor ring (TR) decomposition \cite{zhao2016TensorRing} being the most representative.
TT decomposition can be seen as a special case of TR decomposition, and the latter also possesses other benefits that the former lacks, such as the relaxation of the rank constraint, which results in more intriguing properties, e.g., improved compression capability, better interpretability, and rotational invariance ability. 
So, we mainly focus on TR decomposition in the present paper. For a tensor $\tensor{X} \in \bb{R}^{I_1 \times I_2 \times \cdots \times I_N}$, it can be represented in the element-wise form  as follows:
\begin{align*}
	\tensor{X}(i_1, \cdots, i_N)&=\Trace \left( \mat{G}_1(i_1) \mat{G}_2(i_2) \cdots \mat{G}_N(i_N) \right)
	=\Trace\left( \prod_{n=1}^N \mat{G}_n(i_n) \right),
\end{align*}
where $\mat{G}_n(i_n) =\tensor{G}_n(:,i_n,:) \in \bb{R}^{R_n \times R_{n+1}}$ is the $i_n$-th \emph{lateral slice} of the \emph{core tensor (TR-core)} $\tensor{G}_n \in \bb{R}^{R_n \times I_n \times R_{n+1}}$ with $ R_{N+1}=R_{1}$. Note that a \emph{slice} is a 2nd-order section, i.e., a matrix, of a tensor obtained by fixing all the tensor indices but two. The sizes of TR-cores, i.e., $R_k$ with $k=1,\cdots,N$, are called \emph{TR-ranks}. Additionally, we use the notation $\TR (\{\tensor{G}_n\}_{n=1}^N)$ to denote the TR decomposition of a tensor. 
The problem of fitting $\TR(\{\tensor{G}_n\}_{n=1}^N)$ to a tensor $\tensor{X}$ can be written as the following minimization problem:
\begin{equation}
\label{eq:trmin}
	\mathop{\arg\min}_{\tensor{G}_1, \cdots, \tensor{G}_N} \left\| \TR(\{\tensor{G}_n\}_{n=1}^N) - \tensor{X} \right\|_F,
\end{equation}
where $\| \cdot \|_F$ denotes the Frobenius norm of a matrix or tensor. 
For the above problem, Zhao et al. \cite{zhao2016TensorRing} proposed some algorithms. They can be mainly divided into two categories: one is SVD-based, and the other is the alternating least squares (ALS). 
We will use %the abbreviation %TR-SVD and 
TR-ALS for the latter %method % two methods 
in the following text.

Since both the storage and computational costs increase exponentially as the tensor order grows, the processing difficulties of the aforementioned two approaches will be insurmountable. 
Hence, the development of randomized algorithms for computing TR decompositions of large-scale tensors is becoming increasingly significant; see e.g., \cite{yuan2019RandomizedTensor,ahmadi-asl2020RandomizedAlgorithms,malik2021SamplingBasedMethod,malik2022MoreEfficient,yu2022PracticalSketchingBased}. 
Among them, Malik and Becker \cite{malik2021SamplingBasedMethod} proposed to use the leverage-based sampling to reduce the size of the ALS subproblems and devised a randomized algorithm called TR-ALS-Sampled; Malik \cite{malik2022MoreEfficient} further provided a new approach to approximate the leverage scores and found a more efficient algorithm; 
Yu and Li \cite{yu2022PracticalSketchingBased} proposed two sketching-based randomized algorithms called TR-KSRFT-ALS and TR-TS-ALS.
Most of these algorithms build on the ideas from randomized numerical linear algebra (RandNLA), whose origins lie in theoretical computer science. 

On the other hand, stochastic optimization methods such as the stochastic gradient descent (SGD) have also been employed to develop the algorithms for some tensor decompositions; see e.g., \cite{vervliet2016RandomizedBlock,fu2020BlockRandomizedStochastic,maehara2016ExpectedTensor} for CP decomposition, \cite{li2020SGDTucker} for Tucker decomposition, and \cite{yuan2019HighorderTensor} for TT decomposition.
However, for TR decomposition, there is no related work. To our best knowledge, Yuan et al. \cite{yuan2018HigherdimensionTensor} ever proposed an algorithm called tensor ring weighted optimization (TR-WOPT) for tensor completion problems, which adopts some gradient-based optimization methods;
He and Atia \cite{he2022PatchTrackingbased} later considered the scaled steepest descent method to improve convergence.
However, both of the two methods perform the gradient-based algorithms brutally for the ALS subproblems and hence cannot avoid the formation of oversized coefficient matrices.

In this paper, inspired by \cite{fu2020BlockRandomizedStochastic}, we try to combine the above two perspectives, i.e., the randomized alternating minimization and stochastic gradient-based methods, to propose a doubly randomized optimization framework for TR decomposition, i.e., the block-randomized mini-batch SGD (BRSGD), which can also be regarded as a sensible mix of randomized block coordinate descent (RBCD) and mini-batch SGD. 
To improve the convergence, especially  for ill-conditioned problems, we also consider its scaled version, which borrows the idea of scaled gradient descent (ScaledGD).  It is known that ScaledGD is especially suitable for accelerating the ill-conditioned low-rank matrix estimations, which was analyzed in \cite{tong2021AcceleratingIllconditioned} thoroughly and later was generalized to the low-rank tensor estimation based on Tucker decomposition \cite{tong2021ScalingScalability}. 
To select the right mini-batches of our methods, we propose a theoretical probability distribution with minimum variance, and consider three practical ones, %probability distributions, 
i.e., the uniform, % probability distribution, the
leverage-based, % probability distribution, 
and Euclidean-based probability distributions. 
%It is worth emphasizing that it is not necessary to form the related large matrices when finding these probability distributions. 

The remainder of this paper is organized as follows. 
\Cref{sec:preliminaries} mainly introduces some operations and basics for TR decomposition. The main methods and their convergence analysis, as well as four probability distributions, are presented in \Cref{sec:algorithms}. In \Cref{sec:experiments}, we discuss the numerical performance of our methods. Finally,  the concluding remarks of the whole paper are provided. All proofs are delivered to the Appendix.

%%=================%%
%%  Preliminaries  %%
%%=================%%
\section{Preliminaries}
\label{sec:preliminaries}
We first introduce some definitions on tensor operations, which are necessary later in this paper. 

\begin{definition}[Multi-index]
	\label{def:idx}
	For a positive integer $I$, let $[I] := \{ 1, \cdots, I \}$. 
	For the indices $i_1 \in [I_1], \cdots, i_N \in [I_N]$, a multi-index $i = \overline{i_1 i_2 \cdots i_N}$ refers to an index which takes all possible combinations of values of the indices, $i_1, i_2, \cdots, i_N$, for $i_n = 1,2, \cdots, I_n$ with $n = 1,2, \cdots, N$ in the following %a specific 
    order:
	\begin{equation*}
		\overline{i_1 i_2 \cdots i_N} = i_1+(i_2-1)I_1 +(i_3-1)I_1 I_2+ \cdots +(i_N-1)I_1 \cdots I_{N-1}.
	\end{equation*}  
\end{definition}

\begin{definition}[Mode-$n$ Unfolding]
\label{def:moden}
	The \textbf{mode-$n$ unfolding} of a tensor $\tensor{X} \in \bb{R}^{I_1 \times I_2 \cdots \times I_N}$ is the matrix $\mat{X}_{[n]}$ of size $I_n \times \prod_{k \ne n} I_k$ defined element-wise via
	\begin{equation*}
		\mat{X}_{[n]}(i_n, \overline{i_{n+1} \cdots i_N i_1 \cdots i_{n-1}})=\tensor{X}(i_1, \cdots, i_N).
	\end{equation*}
\end{definition}

\begin{definition}[Classical Mode-$n$ Unfolding]
	The \textbf{classical mode-$n$ unfolding} of a tensor $\tensor{X} \in \bb{R}^{I_1 \times I_2 \cdots \times I_N}$ is the matrix $\mat{X}_{(n)}$ of size $I_n \times \prod_{k \ne n} I_k$ defined element-wise via
	\begin{equation*}
		\mat{X}_{(n)}(i_n, \overline{i_1 \cdots i_{n-1} i_{n+1} \cdots i_N})=\tensor{X}(i_1, \cdots, i_N).
	\end{equation*}
\end{definition}

\begin{definition} [Subchain Product \cite{yu2022PracticalSketchingBased}] 
	\label{def:subchain_product}
	Let $\tensor{A} \in \bb{R}^{I_1 \times J_1 \times K}$ and $\tensor{B} \in \bb{R}^{K \times J_2 \times I_2}$ be two 3rd-order tensors, and $\mat{A}(j_1)$ and $\mat{B}(j_2)$ be the $j_1$-th and $j_2$-th lateral slices of $\tensor{A}$ and $\tensor{B}$, respectively. The mode-2 \textbf{subchain product} of $\tensor{A}$ and $\tensor{B}$ is a tensor of size $I_1 \times J_1 J_2 \times I_2$ denoted by $\tensor{A} \boxtimes_2 \tensor{B}$ and  defined as 
	\begin{equation*}
		(\tensor{A} \boxtimes_2 \tensor{B})(\overline{j_1 j_2}) = \tensor{A}(j_1)\tensor{B}(j_2).
	\end{equation*}
	That is, with respect to the correspondence on indices, the lateral slices of $\tensor{A} \boxtimes_2 \tensor{B}$ are the classical matrix products of the lateral slices of $\tensor{A}$ and $\tensor{B}$. %The mode-1 and mode-3 subchain products can be defined similarly.
\end{definition}

\begin{definition}[Subchain Tensor]
\label{def:subchain}
	Let $\tensor{X} = \TR(\{\tensor{G}_n\}_{n=1}^N)$ be the TR decomposition of $\tensor{X} \in \bb{R}^{I_1 \times I_2 \cdots \times I_N}$. The \textbf{subchain tensor} $\tensor{G}^{\ne n} \in \bb{R}^{R_{n+1} \times \prod_{k \ne n} I_k \times R_n}$ is the merging of all TR-cores expect the $n$-th one and can be written via subchain product as
	\begin{equation*}
		\tensor{G}^{\ne n} = \tensor{G}_{n+1} \boxtimes_2 \cdots \boxtimes_2 \tensor{G}_{N} \boxtimes_2 \tensor{G}_{1} \boxtimes_2 \cdots \boxtimes_2 \tensor{G}_{n-1}.
	\end{equation*}
\end{definition}

\begin{definition} [Slices-Hadamard Product \cite{yu2022PracticalSketchingBased}]
\label{def:hada_product}
	Let $\tensor{A} \in \bb{R}^{I_1 \times J \times K}$ and $\tensor{B} \in \bb{R}^{K \times J \times I_2}$ be two 3rd-order tensors, and $\mat{A}(j)$ and $\mat{B}(j)$ are the $j$-th lateral slices of $\tensor{A}$ and $\tensor{B}$, respectively. The mode-2 \textbf{slices-Hadamard product} of $\tensor{A}$ and $\tensor{B}$ is a tensor of size $I_1 \times J \times I_2$ denoted by $\tensor{A} \boxast_2 \tensor{B}$ and defined as
	\begin{equation*}
		(\tensor{A} \boxast_2 \tensor{B})(j) = \tensor{A}(j)\tensor{B}(j).
	\end{equation*}
	That is, the $j$-th lateral slice of $\tensor{A} \boxast_2 \tensor{B}$ is the classical matrix product of the $j$-th lateral slices of $\tensor{A}$ and $\tensor{B}$. %The mode-1 and mode-3 slices-Hadamard products can be defined similarly.
\end{definition}

Next, we recall two existing methods and the related technique for computing TR decomposition, which are bases for our work.
Firstly, we reformulate the problem \eqref{eq:trmin} %can be reformulated 
using the following optimization model:
\begin{equation}
\label{eq:trmin_opt}
	\mathop{\arg\min}_{\tensor{G}_1, \cdots, \tensor{G}_N} f(\tensor{Y}),
\end{equation}
where $\tensor{Y} := (\tensor{G}_1, \cdots, \tensor{G}_N)$ and $f(\tensor{Y}) = \frac{1}{2} \left\| \TR(\{\tensor{G}_n\}_{n=1}^N) - \tensor{X} \right\|_F^2$.
Further, according to \cite[Theorem 3.5]{zhao2016TensorRing}, the above model can be rewritten as $N$ independent subproblems:
\begin{equation}
\label{eq:tr_als_opt}
\mathop{\arg\min}_{\tensor{G}_n} f(\tensor{Y}),\ n=1,\cdots, N,
\end{equation}
where $f(\tensor{Y}) = \frac{1}{2} \left\| \mat{X}_{[n]} - \mat{G}_{n(2)} (\mat{G}^{\ne n}_{[2]})^\intercal \right\|_F^2$. Thus, ALS can be used to solve \eqref{eq:tr_als_opt} and the specific algorithm is summarized in \Cref{alg:tr_als}.

\begin{algorithm}
	\caption{TR-ALS \cite{zhao2016TensorRing}}
	\label{alg:tr_als}
	\begin{algorithmic}[1]
		\Function{$\{\tensor{G}_n\}_{n=1}^N$= TR-ALS}{$\tensor{X}, R_1, \cdots, R_N$} 
		\Comment $\tensor{X}$ is the input tensor
		
		\Comment $(R_1, \cdots, R_N)$ are the TR-ranks
		\State Initialize cores $\tensor{G}_2, \cdots, \tensor{G}_N$ \label{line:als_init}
		\Repeat
		\For{$n = 1, \cdots, N$}
			\State Compute $\mat{G}_{[2]}^{\ne n}$ from cores \label{line:als_subchain}
			\State Update $\tensor{G}_n = \mathop{\arg\min}_{\tensor{Z}} \left\| \mat{G}_{[2]}^{\ne n} \mat{Z}_{(2)}^\intercal - \mat{X}_{[n]}^\intercal \right\|_F$ \label{line:als_ls}
		\EndFor
		\Until{termination criteria met}
		\State \Return  $\tensor{G}_1, \cdots, \tensor{G}_N$
		\EndFunction
	\end{algorithmic}
\end{algorithm}

Now, assume that one samples a set of rows from $\mat{G}^{\ne n}_{[2]}$ %mode-$n$ fibers 
indexed by $\fn \subset \{1, \cdots, J_n \}$  with $J_n = \prod_{k \ne n}I_k$, which also corresponds to the indices of  mode-$n$ fibers %the rows sampled 
from $\tensor{X}$. %$\mat{G}^{\ne n}_{[2]}$. 
Then, the sampled ALS subproblem of \Cref{line:als_ls} in  \Cref{alg:tr_als} is:
\begin{equation*}
	\tensor{G}_n = \mathop{\arg\min}_{\tensor{Z}} \left\| \mat{G}_{[2]}^{\ne n}(\fn,:) \mat{Z}_{(2)}^\intercal - \mat{X}_{[n]}^\intercal(\fn,:) \right\|_F.
\end{equation*}
With a technique on importance sampling without forming the large matrix $\mat{G}_{[2]}^{\ne n}$, Malik and Becker \cite{malik2021SamplingBasedMethod} proposed a random sampling variant of TR-ALS, i.e., the TR-ALS-Sampled. 
The technique is mainly due to  \Cref{def:subchain} and the corresponding algorithm is summarized in \Cref{alg:sstp-st}.

\begin{algorithm}
\caption{Sampled subchain and data tensor with probabilities 
 (SSDTP), summarized from \cite{malik2021SamplingBasedMethod}}
\label{alg:sstp-st}
\begin{algorithmic}[1]
    \Function{[$\tensor{G}^{\ne n}(:,\fn,:)$, $\mat{X}_{[n]}(:,\fn)$, $\vect{p}_{\fn}$]= SSDTP}{$\{\tensor{G}_k\}_{k=1,k\ne n}^N, |\fn|, \{\vect{p}_k\}_{k=1,k\ne n}^N$} 
    
    \Comment $\tensor{G}_k \in \bb{R}^{R_k \times I_k \times R_{k+1}}$
    
    \Comment $|\fn|$ is the sampling size 
     
    \Comment $\vect{p}_k\in \bb{R}^{I_k}$ is the probability distribution for $\mat{G}_{k(2)}$
    
    \State $\texttt{idxs}=\textsc{Zeros}(|\fn|, N-1)$, $\vect{p}_{\fn} = \textsc{Ones}(|\fn|,1)$
    \For{$k = n+1, \cdots, N, 1, \cdots, n-1$}
        \State $\texttt{idxs}(:,k) = \textsc{Randsample}(I_k, |\fn|, true, \vect{p}_k)$
        \State $\vect{p}_{\fn} = \vect{p}_{\fn} \ast \vect{p}_k(\texttt{idxs}(:,k),:)$ \Comment $\ast$ denotes the Hardmard product
    \EndFor
    \State Set $\tensor{G}^{\ne n}(:,\fn,:) \in \bb{R}^{R_{n+1} \times m \times R_{n+1}}$ %be a tensor of size $R_{n} \times m \times R_{n}$ 
    with each lateral slice being an $R_{n+1} \times R_{n+1}$ identity matrix
    \For{$k = n+1, \cdots, N, 1, \cdots, n-1$}
        \State $\tensor{G}^{\ne n}(:,\fn,:) \leftarrow \tensor{G}^{\ne n}(:,\fn,:) \boxast_2 \tensor{G}_{k}(:, \texttt{idxs}(:,k), :)$
    \EndFor
    \State $\mat{X}_{[n]}(:,\fn) = \textsc{Mode-n-Unfolding}(\tensor{X}(\texttt{idxs}(:,1), \cdots, \texttt{idxs}(:,n-1),:,\texttt{idxs}(:,n+1),\cdots \texttt{idxs}(:,N)))$
    \State \Return  $\tensor{G}^{\ne n}(:,\fn,:)$, $\mat{X}_{[n]}(:,\fn)$ and $\vect{p}_{\fn}$
    \EndFunction
\end{algorithmic}
\end{algorithm}

In addition, the following %two
definition is also necessary throughout the rest of this paper.

\begin{definition}
	We say $\vect{p} \in [0,1]^N$, i.e., an $N$ dimensional vector with the entries being in $[0, 1]$, is a \textbf{probability distribution} if 
	$\sum_{i=1}^N p_i = 1$. 
\end{definition}

%\begin{definition}
	%For a random variable $\xi \in [N]$, we say $\xi \sim \Multi(\vect{p})$ if $\vect{p} \in [0,1]^N$ is a probability distribution on $[N]$ and ${\Pr}(\xi = i) = p_i$.
%\end{definition}

%%=================%%
%%    Algorithms   %%
%%=================%%
\section{Proposed Methods}
\label{sec:algorithms}
We first deduce the partial derivatives of the objective function \eqref{eq:trmin_opt} w.r.t. $\mat{G}_{n(2)}$ as follows:
\begin{equation}
\label{eq:gradient}
	\nabla_{\mat{G}_{n(2)}}f(\tensor{Y}) = \mat{G}_{n(2)} (\mat{G}^{\ne n}_{[2]})^\intercal \mat{G}^{\ne n}_{[2]} - \mat{X}_{[n]} \mat{G}^{\ne n}_{[2]},\ n = 1, \cdots, N.
\end{equation}
Thus, the update rule of  gradient descent algorithm for TR decomposition (TR-GD) can be given as 
\begin{equation*}
	\mat{G}_{n(2)}^{t+1} \leftarrow \mat{G}_{n(2)}^{t} - \alpha^t \nabla_{\mat{G}_{n(2)}}f(\tensor{Y}^t), ~n = 1, \cdots, N,
\end{equation*}
where $\alpha^t$ is the step size. This scheme first appeared in \cite[TR-WOPT]{yuan2018HigherdimensionTensor}.
Here, we consider the mini-batch SGD and combine it with RBCD, i.e., a doubly randomized computational framework, for minimizing \eqref{eq:trmin_opt}. 
In a high-level, at each iteration, we first sample a mode $n$ from all modes of the tensor, and then sample a set of rows of $\mat{G}^{\ne n}_{[2]}$ and the corresponding mode-$n$ fibers of $\tensor{X}$ %indexed by $\fn \subset \{1, \cdots, J_n \}$ 
for the sampled mode %$n$ as done in \Cref{alg:sstp-st} with some probability distribution 
to form an estimate of gradient. Finally, the mini-batch gradient descent is implemented. The details are presented in the following subsections.

%%==================
\subsection{BRSGD for TR decomposition}
\label{ssec:TR-BRSGD}
Assume $\vect{q}^{\ne n} = [p_1, \cdots, p_{J_n}]^\intercal$ is a probability distribution for $\mat{G}^{\ne n}_{[2]}$ and sample $|\fn|$ rows of $\mat{G}^{\ne n}_{[2]}$ %as a batch 
according to $\vect{q}^{\ne n}$. This is can be carried out efficiently by sampling $|\fn|$ slices from each of $\{\tensor{G}_k\}_{k=1,k\ne n}^N$ independently with the probability distribution $\vect{p}_k $ for $\mat{G}_{k(2)}$. % as shown in \Cref{alg:sstp-st}. 
We will revisit this fact in \Cref{ssec:TR-BRSGD-probability1}.  
Then, with $\mat{X}_{[n]}(:,\fn)$, i.e., sampling $|\fn|$ mode-$n$ fibers of $\tensor{X}$ using the newfound index set $\fn$, we can construct an estimate of $\nabla_{\mat{G}_{n(2)}}f(\tensor{Y})$ as follows
\begin{equation}
\label{eq:S_gradient} 
	\mf{g}_n = \frac{1}{|\fn| J_n} \left(\mat{G}_{n(2)} (\mat{G}^{\ne n}_{[2]}(\fn,:))^\intercal \mat{D} \mat{G}^{\ne n}_{[2]}(\fn,:) - \mat{X}_{[n]}(:,\fn) \mat{D} \mat{G}^{\ne n}_{[2]}(\fn,:) \right),
\end{equation}
where $\mat{D} = {\rm diag} \left( \frac{1}{p_{j_1}}, \cdots, \frac{1}{p_{j_{|\fn|}}} \right)$ with $j_1, \cdots, j_{|\fn|} \in \fn$. 
Thus, setting the search direction of the $t$-th iteration $\mf{p}_n^t = -\mf{g}_n^t$, the latent variables can be updated by
\begin{align}
\label{eq:update_brsgd} &\mat{G}_{n(2)}^{t+1} \leftarrow \mat{G}_{n(2)}^{t} + \alpha^t \mf{p}_n^t, ~n = 1, \cdots, N, \\
\nonumber &\mat{G}_{n'(2)}^{t+1} \leftarrow \mat{G}_{n'(2)}^{t} ,~~~ n' \ne n.
\end{align}
The proposed update is very efficient since the most resource-consuming terms $\mat{X}_{[n]} \mat{G}^{\ne n}_{[2]}$ and $(\mat{G}^{\ne n}_{[2]})^\intercal \mat{G}^{\ne n}_{[2]}$ in \eqref{eq:gradient} are avoided. Now, the corresponding parts $\mat{X}_{[n]}(:,\fn) \mat{G}^{\ne n}_{[2]}(\fn,:)$ and $(\mat{G}^{\ne n}_{[2]}(\fn,:))^\intercal \mat{G}^{\ne n}_{[2]}(\fn,:)$ only cost $\bigO{|\fn| I_n R_nR_{n+1}}$ and $\bigO{|\fn| R_n^2 R_{n+1}^2}$ with $|\fn|$ being the input parameter under control. Combining with choosing the mode $n$ uniformly, we summarize the algorithm  framework in \Cref{alg:TR_BRSGD} and call it TR-BRSGD. For \Cref{line:sgd_ssdtp} in \Cref{alg:TR_BRSGD}, we will explain it in \Cref{ssec:TR-BRSGD-probability1}. 

\begin{algorithm} \small
\caption{TR-BRSGD (Proposed Algorithm)}
\label{alg:TR_BRSGD}
\begin{algorithmic}[1]
    \Function{$\{ \tensor{G}_{n} \}_{n=1}^N$= TR-BRSGD}{$\tensor{X}, |\fn|, R_1, \cdots, R_N, \{\alpha^t\}_{t=0,1, \cdots}, \{\vect{p}_k\}_{k=1,k\ne n}^N$} 
    
    \Comment $\tensor{G}_n \in \bb{R}^{R_n \times I_n \times R_{n+1}}$; $\tensor{X} \in \bb{R}^{I_1 \times \cdots \times I_N}$

    \Comment $|\fn|$ is the sampling size 

        \Comment $(R_1, \cdots, R_N)$ are the TR-ranks
    
    \Comment $\{\alpha^t\}_{t=0,1, \cdots}$ are step sizes

    \Comment $\vect{p}_k\in \bb{R}^{I_k}$ is the probability distribution for $\mat{G}_{k(2)}$
    
    \State Initialize cores ${\tensor{G}}_1, \cdots, {\tensor{G}}_N$
    \State $t\leftarrow 0$
    \Repeat
    \State Uniformly sample $n$ from $\{1,\cdots,N\}$
    \State [$\tensor{G}^{\ne n}(:,\fn,:)$, $\mat{X}_{[n]}(:,\fn)$, $\vect{p}_{\fn}$]= \textsc
    {SSDTP}($\{\tensor{G}_k\}_{k=1,k\ne n}^N, |\fn|, \{\vect{p}_k\}_{k=1,k\ne n}^N$) \label{line:sgd_ssdtp}
    \State Compute the search direction $\mf{p}_n^t = -\mf{g}_n^t$ via \eqref{eq:S_gradient} \label{line:sd}
    \State Update $\mat{G}_{n(2)}^{t+1} \leftarrow \mat{G}_{n(2)}^{t} + \alpha^t \mf{p}_n^t$, $\mat{G}_{n'(2)}^{t+1} \leftarrow \mat{G}_{n'(2)}^{t}$ for $n' \ne n$
    \State $t\leftarrow t+1$
    \Until some stopping criterion is reached
    \State \Return $\{ \tensor{G}_{n} \}_{n=1}^N$
    \EndFunction
\end{algorithmic}
\end{algorithm}

\begin{remark}
\label{rem:adagrad}
An adaptive step size scheme that incorporates the AdaGrad can be combined with TR-BRSGD:
\begin{align*}
% \label{eq:update_eta} 
&[\vect{\eta}_n^t]_{i, r} \leftarrow \frac{\eta}{\left(b + \sum_{t'=1}^t[-\mf{p}_n^{t'}]_{i, r}^2\right)^{(1/2+\epsilon)}}, ~~~ i \in [I_n], r \in [R_{n}R_{n+1}], \\
% \label{eq:update_brsgdada}
&\mat{G}_{n(2)}^{t+1} \leftarrow \mat{G}_{n(2)}^{t} + \vect{\eta}_n^t \ast \mf{p}_n^t, ~~~n = 1, \cdots, N, \\
\nonumber &\mat{G}_{n'(2)}^{t+1} \leftarrow \mat{G}_{n'(2)}^{t} ,~~~ n' \ne n,
\end{align*}
where %$\ast$ denotes the Hadamard product, and 
$\eta,b,\epsilon>0$. 
Here, $b, %>0$ and $
\epsilon>0$ are technical conditions. 
In practice, setting $b=\epsilon=0$ does not hurt the performance. The above update rule can tun the step size throughout the current search direction during the running process.
\end{remark}

Next, we consider the convergence of TR-BRSGD. To this end, we rewrite \eqref{eq:tr_als_opt} as
\begin{equation*}
	\mathop{\arg\min}_{\tensor{G}_n} \frac{1}{J_n} \sum_{i=1}^{J_n} f_i(\tensor{Y}),
\end{equation*}
where $f_i(\tensor{Y}) = \frac{J_n}{2} \left\| \mat{X}_{[n]}(:,i) - \mat{G}_{n(2)} (\mat{G}^{\ne n}_{[2]}(i,:))^\intercal \right\|_F^2$, and define $\xi^{t} \in \{1,\cdots,N\}$ and $\zeta^{t} \subseteq \{1,\cdots,\prod_{m \ne \xi^{t}}I_m \}$ as the random variables (r.v.s) responsible for selecting the mode and fibers in the $t$-th iteration, respectively.  Thus, we first have the unbiasedness of the stochastic gradient given in \eqref{eq:S_gradient}.

\begin{lemma}[Unbiased Gradient]
\label{lem:unbias_gradient}
	Denote $\mc{B}^{t}$ as the filtration generated by the r.v.s 
	\begin{equation*}
		\{ \xi^{0},\zeta^{0}, \xi^{1},\zeta^{1}, \cdots,\xi^{t-1},\zeta^{t-1} \}
	\end{equation*}
	such that the $t$-th iteration $\tensor{Y}^{t}$ is determined conditioned on $\mc{B}^{t}$. Then the stochastic gradient $\mf{g}_n^t$  in the $t$-th iteration formed as in \eqref{eq:S_gradient} is the unbiased estimate of the full gradient w.r.t. $\mat{G}_{\xi^{t}(2)}$:
	\begin{equation*}
	\bb{E}_{\zeta^{t}}\left[ \mf{g}_{\xi^{t}}^t~|~\mc{B}^{t},\xi^{t} \right] = \nabla_{\mat{G}_{\xi^{t}(2)}}f(\tensor{Y}^{t}).
	\end{equation*}
\end{lemma}

Based on \Cref{lem:unbias_gradient}, %the unbiasedness of the stochastic gradient given in \eqref{eq:S_gradient}, 
we then have the following convergence property for TR-BRSGD. % algorithm. %whose proof can be found in \Cref{appdxpf:thm:conv_brsgd}.

\begin{theorem}
\label{thm:conv_brsgd}
Suppose that the step size schedule follows the Robbins-Monro rule: $\sum_{t=0}^{\infty} \alpha^t = \infty$ and $\sum_{t=0}^{\infty} (\alpha^t)^2 < \infty$, and the updates $\mat{G}_{n(2)}^{t}$ are bounded for all $n$ and $t$. Then the sequence produced by TR-BRSGD satisfies:
\begin{equation*}
	\liminf_{t \rightarrow \infty} \bb{E} \left[ \left\| \nabla f(\tensor{Y}^{t}) \right\|_F^2 \right] = 0.
\end{equation*}
\end{theorem}

%%==================
\subsection{Four probability distributions}
\label{ssec:TR-BRSGD-probability}
Considering that the probability distribution $\vect{q}^{\ne n}$ or $\vect{p}_k $ plays an important role in BRSGD, in this subsection, we will present a probability distribution with minimum variance called optimal probability distribution and also discuss three practical ones.

\subsubsection{Optimal probability distribution}
The optimal probability distribution can make the variance of the gradient estimation error $\delta_{\xi^{t}}^t =  \mf{g}_{\xi^{t}}^t - \nabla_{\mat{G}_{\xi^{t}(2)}}f(\tensor{Y}^{t})$ be minimum, which in turn can speed up the convergence rate. 
The results are summarized in \Cref{thm:variance}. % whose proof can be found in \Cref{appdxpf:variance}.

\begin{theorem}
\label{thm:variance}
Denote $\mc{B}^{t}$ as the filtration generated by the r.v.s 
$$\{ \xi^{0},\zeta^{0}, \xi^{1},\zeta^{1}, \cdots,\xi^{t-1},\zeta^{t-1} \}$$ 
such that the $t$-th iteration $\tensor{Y}^{t}$ is determined conditioned on $\mc{B}^{t}$. Suppose that $\vect{q}^{\ne {\xi^{t}}} = [p_1, \cdots, p_{J_{\xi^{t}}}]^\intercal$ is a probability distribution and $\mf{g}_{\xi^{t}}^t$ is computed as in \eqref{eq:S_gradient}. Then,
{\footnotesize\begin{align}
\label{eq:var}
	\bb{E}_{\zeta^{t}}\left[ \left\| \mf{g}_{\xi^{t}}^t - \nabla_{\mat{G}_{\xi^{t}(2)}}f(\tensor{Y}^{t}) \right\|^2_F~|~\mc{B}^{t},\xi^{t} \right] 
	&\nonumber = \frac{1}{|\fxi|} \sum_{j_f = 1}^{J_{\xi^{t}}} \frac{1}{p_{j_f}} \left\|\mat{R}_{[\xi^{t}]}^t(:,j_f) \right\|_2^2 \left\|(\mat{G}^{\ne \xi^{t}}_{[2]})^t(j_f,:)\right\|_2^2 \\
	&~~~~ - \frac{1}{|\fxi|} \left\| \nabla_{\mat{G}_{\xi^{t}(2)}}f(\tensor{G}_{\xi^{t}}) \right\|_F^2.
\end{align}}
Further, if $\vect{q}^{\ne \xi^{t}} $ is as
\begin{equation}
\label{eq:optimal_prob}
	\vect{q}^{\ne \xi^{t}}(i) = \frac{\left\| \mat{R}_{[\xi^{t}]}^t(:,i) \right\|_2 \left\|(\mat{G}^{\ne \xi^{t}}_{[2]})^t(i,:)\right\|_2}{\sum_{i'=1}^{J_{\xi^{t}}} \left\|\mat{R}_{[\xi^{t}]}^t(:,i) \right\|_2 \left\|(\mat{G}^{\ne \xi^{t}}_{[2]})^t(i,:)\right\|_2},
\end{equation}
where $\mat{R}_{[\xi^{t}]}^t = \mat{X}_{[\xi^{t}]}^t - \mat{X}_{[\xi^{t}]}$, $\mat{X}_{[\xi^{t}]}^t = \mat{G}_{\xi^{t}(2)}^t ((\mat{G}^{\ne \xi^{t}}_{[2]})^t)^\intercal$, and $i \in [J_{\xi^{t}}]$,
then {\footnotesize$\bb{E}_{\zeta^{t}}\left[ \left\| \mf{g}_{\xi^{t}}^t - \nabla_{\mat{G}_{\xi^{t}(2)}}f(\tensor{Y}^{t}) \right\|^2_F~|~\mc{B}^{t},\xi^{t} \right]$} achieves its minimum as
{\footnotesize\begin{align}
\label{eq:optimal_var}
	\bb{E}_{\zeta^{t}}\left[ \left\| \mf{g}_{\xi^{t}}^t - \nabla_{\mat{G}_{\xi^{t}(2)}}f(\tensor{Y}^{t}) \right\|^2_F~|~\mc{B}^{t},\xi^{t} \right]
	&\nonumber= \frac{1}{|\fxi|} \left( \sum_{j_f = 1}^{J_{\xi^{t}}} \left\|\mat{R}_{[\xi^{t}]}^t(:,j_f) \right\|_2 \left\|(\mat{G}^{\ne \xi^{t}}_{[2]})^t(j_f,:)\right\|_2 \right)^2 \\
	&~~~~ - \frac{1}{|\fxi|} \left\| \nabla_{\mat{G}_{\xi^{t}(2)}}f(\tensor{G}_{\xi^{t}}) \right\|_F^2.
\end{align}}
\end{theorem}

Unfortunately, the optimal probabilities in \eqref{eq:optimal_prob} are unpractical since 
the full matrix $\mat{R}_{[\xi^{t}]}^t$ needs to be formed to compute these probabilities.  
In the following subsections, we consider three practical probability distributions, which can also be combined with the sampling method in \Cref{alg:sstp-st}.

\subsubsection{Practical probability distributions}
\label{ssec:TR-BRSGD-probability1}
We first consider the \emph{uniform probability distribution.} That is, 
\begin{equation*}
	\vect{q}^{\ne n} = \left[\frac{1}{J_n}, \cdots, \frac{1}{J_n}\right]^\intercal \in \bb{R}^{J_n}.
\end{equation*}
In this case, the expression \eqref{eq:S_gradient} can be simplified as follows
\begin{equation*}
	\mf{g}_n = \frac{1}{|\fn|} \left(\mat{G}_{n(2)} (\mat{G}^{\ne n}_{[2]}(\fn,:))^\intercal \mat{G}^{\ne n}_{[2]}(\fn,:) - \mat{X}_{[n]}(:,\fn) \mat{G}^{\ne n}_{[2]}(\fn,:) \right).
\end{equation*}

Next, we investigate the \emph{Leverage-based probability distribution.}
Two definitions are first introduced.
\begin{definition}[Leverage Scores \cite{drineas2012FastApproximation}]
\label{def:leverages_scores}
	Let $\mat{A}\in \bb{R}^{m\times n}$ with $m > n$, and
	let $\mat{Q} \in \bb{R}^{m \times n}$ be any orthogonal basis for the column space of $\mat{A}$.
	The \textbf{leverage score} of the $i$-th row of $\mat{A}$ is given by
	\begin{equation*}
		\ell_i(\mat{A}) = \|\mat{Q}(i,:)\|_2^2.
	\end{equation*}
\end{definition}

\begin{definition}[Leverage-based Probability Distribution \cite{woodruff2014SketchingTool}]
\label{def:leverages_scores_sampling}
	Let $\mat{A}\in \bb{R}^{m\times n}$ with $m > n$.
	We say a probability distribution $\vect{q}=[q_1, \cdots, q_{m} ]^\intercal$ is a \textbf{leverage-based probability distribution} for $\mat{A}$ on $[m]$ if $q_i \ge \beta p_i$ with $p_i=\frac{\ell_i(\mat{A})}{n}$, $0<\beta \le 1$ and $\forall i \in [m]$.
\end{definition}

Since it is expensive to compute the leverage scores of $\mat{G}^{\ne n}_{[2]} \in \bb{R}^{J_n \times R_{n}R_{n+1}}$ directly, which costs $\bigO{J_n R_{n}^2R_{n+1}^2}$, according to \cite{malik2021SamplingBasedMethod}, we can estimate them from the leverage scores related to the cores $\tensor{G}_{1}, \cdots, \tensor{G}_{n-1},\tensor{G}_{n+1}, \cdots, \tensor{G}_{N}$.

\begin{lemma}[\cite{malik2021SamplingBasedMethod}]
\label{lem:lev_est}
	For each $n \in [N]$, let $\vect{p}_n \in \bb{R}^{I_n}$ be a probability distribution on $[I_n]$ defined element-wise via
	\begin{equation*}
		\vect{p}_n(i_n) = \frac{\ell_{i_n}(\mat{G}_{n(2)})}{\rank(\mat{G}_{n(2)}))},
	\end{equation*}
	$\vect{p}^{\ne n}\in \bb{R}^{J_n}$ be a probability distribution on $\left[J_n\right]$ defined element-wise via
	\begin{equation*}
		\vect{p}^{\ne n}(i) = \frac{\ell_{i}(\mat{G}^{\ne n}_{[2]})}{\rank(\mat{G}^{\ne n}_{[2]})},
	\end{equation*}
	$\vect{q}^{\ne n}\in \bb{R}^{J_n}$ be a vector defined element-wise via
	\begin{equation*}
		\vect{q}^{\ne n}(\overline{i_{n+1} \cdots i_{N} i_{1} \cdots i_{n-1}}) = \prod_{\substack{k=1\\ k \ne n}}^{N} \vect{p}_k(i_k),
	\end{equation*} 
	and $\beta_n$ be a constant as in \Cref{def:leverages_scores_sampling} defined as 
        % {\footnotesize
	\begin{equation*}
		\beta_{n} = \frac{1}{\left( R_{n}R_{n+1} \prod_{k \notin \{n, n+1\}} R_{k}^2 \right)}. % \substack
	\end{equation*}
        % }
    Then for each $n \in [N]$, $\vect{q}^{\ne n}(i) \geq \beta_n \vect{p}^{\ne n}(i)$ for all $i = \overline{i_{n+1} \cdots i_{N} i_{1} \cdots i_{n-1}} \in \left[J_n\right]$ and hence $\vect{q}^{\ne n}$ is the leverage-based probability distribution for $\mat{G}^{\ne n}_{[2]}$ on $\left[J_n\right]$.
\end{lemma}

Now, we consider the \emph{Euclidean-based probability distribution.}  A definition is first introduced.
\begin{definition}[Euclidean-based Probability Distribution]
\label{def:euclidean_norms_sampling}
Let $\tensor{A} \in \bb{R}^{R_{1} \times I \times R_{2}}$. %with $I > R_{1}R_{2}$.
We say  a probability distribution $\vect{p}=[p_1, \cdots, p_{I} ]^\intercal$ is an mode-2 \textbf{Euclidean-based probability distribution }for $\tensor{A}$ on $[I]$ if $p_i \ge \beta ||\tensor{A}(:, i, :)||_F^2 / \| \tensor{A} \|_F^2$ with $0 < \beta \leq 1$ and $\forall i \in [I]$. 
\end{definition}

Similar to \Cref{lem:lev_est}, we also have that the Euclidean-based probability distribution for $\mat{G}^{\ne n}_{[2]}$ can be estimated from the norms related to the cores $\tensor{G}_{1}, \cdots, \tensor{G}_{n-1},\tensor{G}_{n+1}, \cdots, \tensor{G}_{N}$. % see \Cref{appdxpf:euc_est} for proof.
\begin{lemma}
\label{lem:euc_est}
	For each $n \in [N]$, let $\vect{p}_n \in \bb{R}^{I_n}$ be a probability distribution on $[I_n]$ defined element-wise via
	\begin{equation*}
		\vect{p}_n(i_n) = \frac{\| \tensor{G}_{n}(:,i_n,:) \|_F^2}{\| \tensor{G}_{n} \|_F^2},
	\end{equation*}
	$\vect{p}^{\ne n}\in \bb{R}^{J_n}$ be a probability distribution on $\left[J_n\right]$ defined element-wise via
	\begin{equation*}
		\vect{p}^{\ne n}(i) = \frac{\| \tensor{G}^{\ne n}(:,i,:) \|_F^2}{\| \tensor{G}^{\ne n} \|_F^2},
	\end{equation*}
	$\vect{q}^{\ne n}\in \bb{R}^{J_n}$ be a vector defined element-wise via
	\begin{equation*}
		\vect{q}^{\ne n}(\overline{i_{n+1} \cdots i_{N} i_{1} \cdots i_{n-1}}) = \prod_{\substack{k=1\\ k \ne n}}^{N} \vect{p}_k(i_k),
	\end{equation*}
	and $\beta_n$ be a constant as in \Cref{def:leverages_scores_sampling} defined as
	\begin{equation*}
		\beta_n = \frac{1}{\prod_{k \ne n} \| \tensor{G}_{k} \|_F^2}.
	\end{equation*}
	Then for each $n \in [N]$, $\vect{q}^{\ne n}(i) \geq \beta_n \vect{p}^{\ne n}(i)$ for all $i = \overline{i_{n+1} \cdots i_{N} i_{1} \cdots i_{n-1}} \in \left[J_n\right]$ and hence $\vect{q}^{\ne n}$ is the Euclidean-based probability distribution for $\mat{G}^{\ne n}_{[2]}$ on $\left[J_n\right]$.
\end{lemma}

From the fact that $\frac{1}{J_n}=\prod_{k \ne n} \frac{1}{I_k}$, and \Cref{lem:lev_est,lem:euc_est},  and considering \Cref{def:subchain}, we can find that sampling $|\fn|$ rows of $\mat{G}^{\ne n}_{[2]}$ according to $\vect{q}^{\ne n} = [p_1, \cdots, p_{J_n}]^\intercal$ can be indeed carried out by sampling $|\fn|$ slices from each of $\{\tensor{G}_k\}_{k=1,k \ne n}^N$ independently with the %corresponding
probability distribution $\vect{p}_k $ for $\mat{G}_{k(2)}$. % $\vect{p}_n$ on $[I_n]$. 
The specific process can be seen in \Cref{alg:sstp-st}, where we also show how to form $\mat{X}_{[n]}(:,\fn)$ and $\vect{p}_{\fn}$ efficiently. This constitutes \Cref{line:sgd_ssdtp} in \Cref{alg:TR_BRSGD}. 

%%==================
\subsection{Scaled BRSGD for TR decomposition}
\label{ssec:TR-BRSGD-Scaled}
As done in \cite{tong2021AcceleratingIllconditioned,tong2021ScalingScalability}, for the problem \eqref{eq:trmin_opt}, %or\eqref{eq:tr_als_opt},
 the update rule of ScaledGD can be given as
\begin{equation}
\label{eq:scaled_gd}
	\mat{G}_{n(2)}^{t+1} \leftarrow \mat{G}_{n(2)}^{t} - \alpha^t \nabla_{\mat{G}_{n(2)}}f(\tensor{Y}^t)\left( ((\mat{G}^{\ne n}_{[2]})^t)^\intercal (\mat{G}^{\ne n}_{[2]})^t\right)^{-1}, ~n = 1, \cdots, N,
\end{equation}
which is equivalent to 
\begin{equation*}
	{\rm vec}(\tensor{Y}^{t+1}) = {\rm vec}(\tensor{Y}^{t}) - \alpha^t (\mat{H}^{t})^{-1} \nabla_{{\rm vec}(\tensor{Y})}f(\tensor{Y}^t),
\end{equation*}
where $\mat{H}^{t} = {\rm diag} \left( ((\mat{G}^{\ne 1}_{[2]})^t)^\intercal (\mat{G}^{\ne 1}_{[2]})^t \otimes \mat{I}_{I_1}, \cdots, ((\mat{G}^{\ne N}_{[2]})^t)^\intercal (\mat{G}^{\ne N}_{[2]})^t \otimes \mat{I}_{I_N} \right)$. 
Note that the Hessian of the loss function \eqref{eq:tr_als_opt} w.r.t $\mat{G}_{n(2)}$ %\eqref{eq:trmin_opt} 
is
\begin{equation}
\label{eq:hessian}
	 \nabla^2_{{\rm vec}(\mat{G}_{n(2)}),{\rm vec}(\mat{G}_{n(2)})}f(\tensor{Y}^t) = ((\mat{G}^{\ne n}_{[2]})^t)^\intercal (\mat{G}^{\ne n}_{[2]})^t \otimes \mat{I}_{I_n},
\end{equation}
which is the $n$-th diagonal block of the Hessian of the problem \eqref{eq:trmin_opt}. Therefore, as pointed out in \cite{tong2021AcceleratingIllconditioned,tong2021ScalingScalability}, ScaledGD can be regarded as a quasi-Newton method where the pre-conditioner is designed as the inverse of the diagonal approximation of the Hessian.

Next, we consider the stochastic version of \eqref{eq:scaled_gd}. To this end, we first present an estimate of the Hessian in \eqref{eq:hessian} as follows
\begin{equation}
\label{eq:S_hessian}
	\tilde{\mf{h}}_n = \frac{1}{|\fn| J_n} \left((\mat{G}^{\ne n}_{[2]}(\hn,:))^\intercal \mat{D} \mat{G}^{\ne n}_{[2]}(\hn,:) \right) \otimes \mat{I}_{I_n},
\end{equation}
where $\mat{D}$ is similar to the one in \eqref{eq:S_gradient} and $\hn \subset \{1, \cdots, J_n \}$ contains the indices of the rows sampled from $\mat{G}^{\ne n}_{[2]}$. Additionally, $\hn$ can be the same as $\fn$. 
Then, the pre-conditioner in \eqref{eq:scaled_gd} can be designed as
\begin{equation*}
	(\mf{h}_n)^{-1} =\left[ \frac{1}{|\fn| J_n} \left((\mat{G}^{\ne n}_{[2]}(\hn,:))^\intercal \mat{D} \mat{G}^{\ne n}_{[2]}(\hn,:) \right) \right]^{-1}.
\end{equation*}
Thus, letting the search direction of the $t$-th iteration 
\begin{equation}
\label{eq:scaledsg_direction}
\mf{p}_n^t = -\mf{g}_n^t (\mf{h}_n^t)^{-1}
\end{equation}
and then the latent variables be updated by \eqref{eq:update_brsgd}, we have TR-ScaledBRSGD by using the algorithm framework of \Cref{alg:TR_BRSGD}. That is, replace \Cref{line:sd} in \Cref{alg:TR_BRSGD} with $\mf{p}_n^t = -\mf{g}_n^t (\mf{h}_n^t)^{-1}$. Furthermore, we can also use the sampling method in \Cref{alg:sstp-st} to form $\mf{h}_n$, and combine TR-ScaledBRSGD with the adaptive step size as done in \Cref{rem:adagrad}.

\begin{remark}
Considering that the true rank is often unknown in practice, the target rank of the model is usually overspecified. However, this over-parameterized regime may significantly slow down the convergence of local search algorithms. To tackle this issue,  \cite{zhang2021PreconditionedGradient} proposed an inexpensive amended pre-conditioner for ScaledGD. For our TR-ScaledSGD, the amended pre-conditioner is given as
\begin{equation*}
(\mf{h}_n)^{-1} =\left[ \frac{1}{|\fn| J_n} \left((\mat{G}^{\ne n}_{[2]}(\hn,:))^\intercal \mat{D} \mat{G}^{\ne n}_{[2]}(\hn,:) \right) + \eta^t \mat{I} \right]^{-1},
\end{equation*}
where $\eta^t \geq 0$ is the damping parameter specific to the $t$-th iteration, which can also be viewed as a parameter that allows us to interpolate between TR-ScaledSGD (with $\eta^t = 0$) and TR-BRSGD (in the limit $\eta^t \rightarrow \infty$). In addition, the above amendment can also make TR-ScaledSGD avoid inverting near-singular matrices in iterations. %In addition, $\eta^t$ can be viewed as a parameter that allows us to interpolate between TR-ScaledSGD (with $\eta^t = 0$) and TR-BRSGD (in the limit $\eta^t \rightarrow \infty$).
\end{remark}

Similar to the discussion of the unbiasedness of the stochastic gradient given in \eqref{eq:S_gradient}, we now establish the unbiasedness of the stochastic Hessian given in \eqref{eq:S_hessian} and the search direction given in \eqref{eq:scaledsg_direction}. % as follows. % see \Cref{appdxpf:unbias_hessian,appdxpf:unbias_sd} for proofs, respectively.

\begin{lemma}[Unbiased Hessian]
\label{lem:unbias_hessian}
    Denote $\mc{B}^{t}$ as the filtration generated by the r.v.s 
	\begin{equation*}
		\{ \xi^{0},\zeta^{0}, \xi^{1},\zeta^{1}, \cdots,\xi^{t-1},\zeta^{t-1} \}
	\end{equation*}
    such that the $t$-th iteration $\tensor{Y}^{t}$ is determined conditioned on $\mc{B}^{t}$. 
    Then the stochastic Hessian $\tilde{\mf{h}}_{\xi^{t}}$ formed as in \eqref{eq:S_hessian} is the unbiased estimate of the full Hessian w.r.t ${\rm vec}(\mat{G}_{\xi^{t}(2)})$ formed as in \eqref{eq:hessian}:
    \begin{equation*}
    \bb{E}_{\zeta^{t}}\left[ \tilde{\mf{h}}_{\xi^{t}}^t~|~\mc{B}^{t},\xi^{t} \right] = \nabla^2_{{\rm vec}(\mat{G}_{\xi^{t}(2)}),{\rm vec}(\mat{G}_{\xi^{t}(2)})}f(\tensor{Y}^t),
    \end{equation*}
    and also
    \begin{equation*}
    \bb{E}_{\zeta^{t}}\left[ \mf{h}_{\xi^{t}}^t~|~\mc{B}^{t},\xi^{t} \right] = (\mat{G}^{\ne \xi^{t}}_{[2]})^\intercal \mat{G}^{\ne \xi^{t}}_{[2]}.
    \end{equation*}
\end{lemma}

\begin{lemma}[Approximate Unbiaseed Search Direction]
\label{lem:unbias_sd}
    Denote $\mc{B}^{t}$ as the filtration generated by the r.v.s 
	\begin{equation*}
		\{ \xi^{0},\zeta^{0}, \xi^{1},\zeta^{1}, \cdots,\xi^{t-1},\zeta^{t-1} \}
	\end{equation*}
	such that the $t$-th iteration $\tensor{Y}^{t}$ is determined conditioned on $\mc{B}^{t}$, and
    assume that $\nabla^2_{{\rm vec}(\mat{G}_{\xi^{t}(2)}),{\rm vec}(\mat{G}_{\xi^{t}(2)})}f(\tensor{Y}^t) \preceq \mat{I}$. Then the stochastic search direction $\mf{p}_n^t = -\mf{g}_n^t (\mf{h}_n^t)^{-1}$ is approximately unbiased:
	\begin{equation*}
	\bb{E}_{\zeta^{t}}\left[ \mf{p}_{\xi^{t}}^t~|~\mc{B}^{t},\xi^{t} \right] \approx -\nabla_{\mat{G}_{\xi^{t}(2)}}f(\tensor{Y}^t) \left( (\mat{G}^{\ne \xi^{t}}_{[2]})^\intercal \mat{G}^{\ne \xi^{t}}_{[2]} \right)^{-1}.
	\end{equation*}
\end{lemma}
\begin{remark}
The assumption $\nabla^2_{{\rm vec}(\mat{G}_{\xi^{t}(2)}),{\rm vec}(\mat{G}_{\xi^{t}(2)})}f(\tensor{Y}^t) \preceq \mat{I}$ is without loss of generality due to scaling, and is often used in the literature for simplifying the analysis; see e.g., \cite{agarwal2017SecondOrderStochastic}. %which is intended to simplify the proof; 
Furthermore, the result really used in the analysis is actually the spectral radius of Hessian being smaller than 1. So this weaker condition can replace the aforementioned assumption. 
%in the course of the proof this is simply to make its spectral parametrization less than 1 without such a strong assumption, we can also only assume the Hessian is upper bounded by 1.
\end{remark}

Based on the above approximate unbiasedness of the stochastic search direction $\mf{p}_n^t$, we can obtain the convergence of TR-ScaledBRSGD. %and give the proof in \Cref{appdxpf:conv_scaled}.

\begin{theorem}
\label{thm:conv_scaled}
Suppose that the stepsize schedule follows the Robbins-Monro rule: $\sum_{t=0}^{\infty} \alpha^t = \infty$ and $\sum_{t=0}^{\infty} (\alpha^t)^2 < \infty$, the updates $\mat{G}_{n(2)}^{t}$ are bounded for all $n$ and $t$, and $\nabla^2_{{\rm vec}(\mat{G}_{\xi^{t}(2)}),{\rm vec}(\mat{G}_{\xi^{t}(2)})}f(\tensor{Y}^t) \preceq \mat{I}$.
Then the solution sequence produced by TR-ScaledBRSGD satisfies:
\begin{equation*}
	\liminf_{t \rightarrow \infty} \bb{E} \left[ \left\| (\mat{H}^{t})^{-1/2} \nabla f(\tensor{Y}^{t}) \right\|_F^2 \right] = 0.
\end{equation*}
\end{theorem}

%%=================%%
%%   Experiments   %%
%%=================%%
\section{Numerical Experiments}
\label{sec:experiments}
We will first use synthetic data on well-conditioned and ill-conditioned tensors and then employ real-world problems to verify the superior performance of our approaches, particularly in computational cost.
All experiments are run by using MATLAB (Version 2022a) on a computer with an Intel Xeon W-2255 3.7 GHz CPU and 256 GB memory, and we also use the MATLAB Tensor Toolbox \cite{kolda2006TensorToolbox}.
Additionally, all numerical results and plotted quantities are averages over 10 runs.

In our experiments, the following baselines are %have been 
chosen as benchmarks to evaluate the performance.
\begin{itemize}
\item Traditional methods:
\begin{itemize}
    \item TR-ALS \cite{zhao2016TensorRing}: A classical algorithm. % for TR decomposition.
    \item TR-GD \cite{yuan2018HigherdimensionTensor}: A slight variant of TR-WOPT. % for TR decomposition.
    \item TR-ScaledGD \cite{he2022PatchTrackingbased}: The scaled steepest descent method. % for TR decomposition.
\end{itemize}
\item RandNLA methods:
\begin{itemize}
    \item TR-ALS-Sampled \cite{malik2021SamplingBasedMethod}: A random sampling algorithm based on leverage scores.
    \item TR-KSRFT-ALS \cite{yu2022PracticalSketchingBased}: A random projection algorithm based on Kronecker sub-sampled randomized Fourier transform.
    \item TR-TS-ALS \cite{yu2022PracticalSketchingBased}: A random projection algorithm based on TensorSketch.
\end{itemize}
\end{itemize}
In all the methods, either random initialization or spectral initialization can be adopted. In the specific experiments below, we use random initialization for Experiment-\uppercase\expandafter{\romannumeral1} and spectral initialization for others.
Meanwhile,  we employ the relative square error (RSE) defined by
\begin{equation*}
	\text{RSE} = \frac{\left\| \TR \left( \{\hat{\tensor{G}}_n\}_{n=1}^N \right) - \tensor{X}_{true} \right\|_F}{\left\| \tensor{X}_{true} \right\|_F},
\end{equation*}
where %the TR-cores
 $\hat{\tensor{G}}_n$ is computed by various algorithms,  to measure the performance in most experiments.
For image data, peak signal-to-noise (PSNR) is obtained by:
\begin{equation*}
	\text{PSNR} = 10 \log_{10}(255^2/\text{MSE}),
\end{equation*}
where MSE is deduced by:
\begin{equation*}
	\text{MSE} = \frac{\left\| \TR \left( \{\hat{\tensor{G}}_n\}_{n=1}^N \right) - \tensor{X}_{true} \right\|_F}{\textsc{Num}(\tensor{X}_{true})},
\end{equation*}
with $\textsc{Num}(\cdot)$ denoting the number of elements of the tensor.

Furthermore, we implement three stopping  criteria for all the algorithms: the maximum iterations  $T$, the maximum total running time $MT$, and the minimum
 RSE $tol$. When one of the stopping criteria is met, the algorithm will be terminated.
A final point to be declared is the choice of parameters, e.g., the sketching size $J$ for sketching-based algorithms, and the step size $\alpha$ and batch sizes $|\fn|$ and $|\hn|$ for gradient-based algorithms. 
For the sake of fairness, we select the suggested parameter settings of each algorithm; see the references of these algorithms for guidelines for the selection of parameters.

%%==================
\subsection{Synthetic data}

All the synthetic tensors have the same dimensions in all modes. They are generated by creating $N$ TR-cores of size $R_{true} \times I \times R_{true}$ first, and then formed by $\tensor{X}_{true} = \TR \left( \{\tensor{G}_n\}_{n=1}^N \right)$. 
Note that the TR-cores may be generated in different ways, which will be detailed in subsequent experiments. In addition, we denote the target rank of all algorithms by $R$.

\begin{figure}[htbp] 
	\centering 
	\subfloat[Number of iterations v.s. RSE]{\includegraphics[scale=0.3]{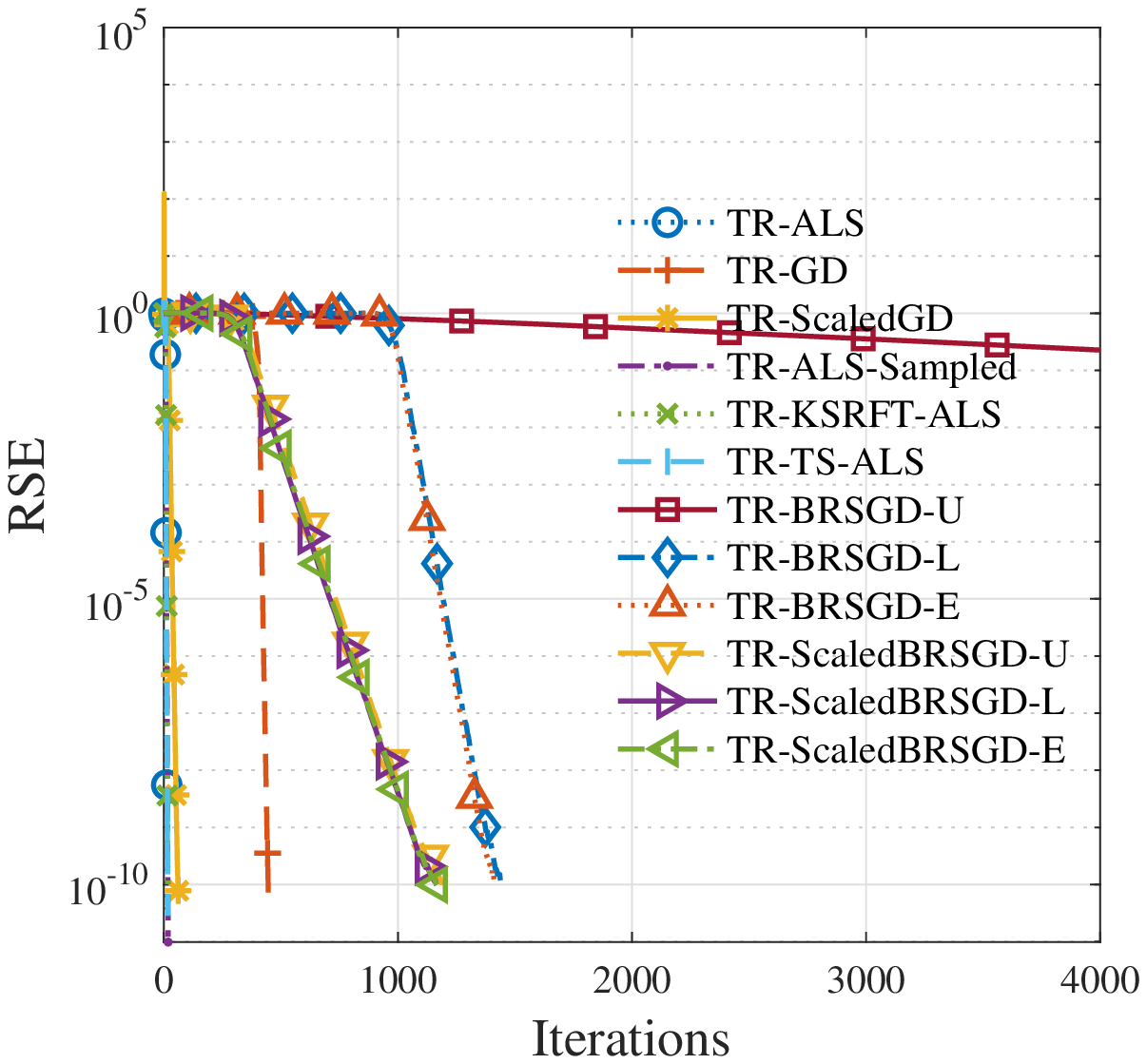}} 
	\subfloat[Time v.s. RSE]{\includegraphics[scale=0.3]{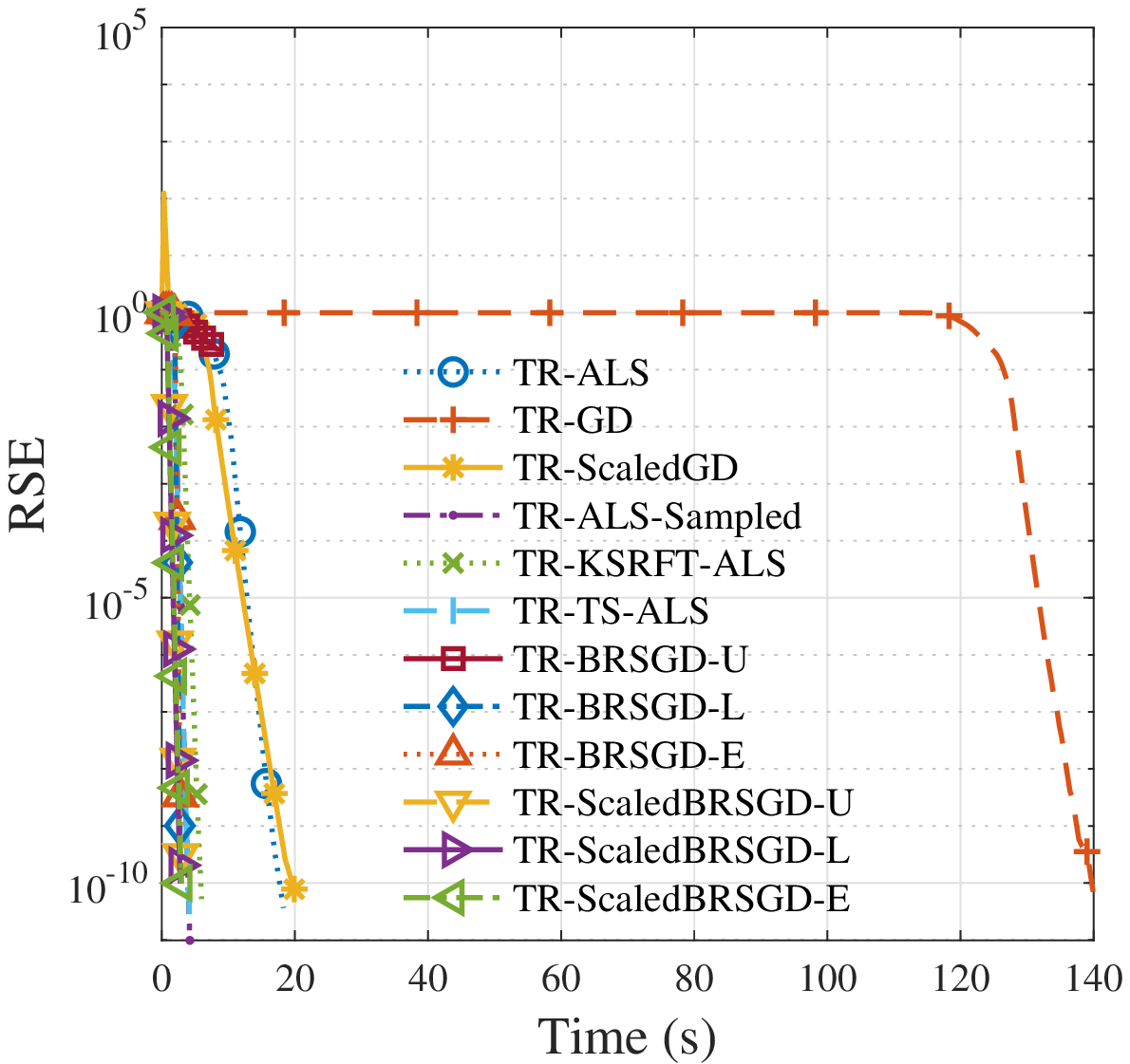}} 
	\subfloat[Zoom in Time v.s. RSE]{\includegraphics[scale=0.3]{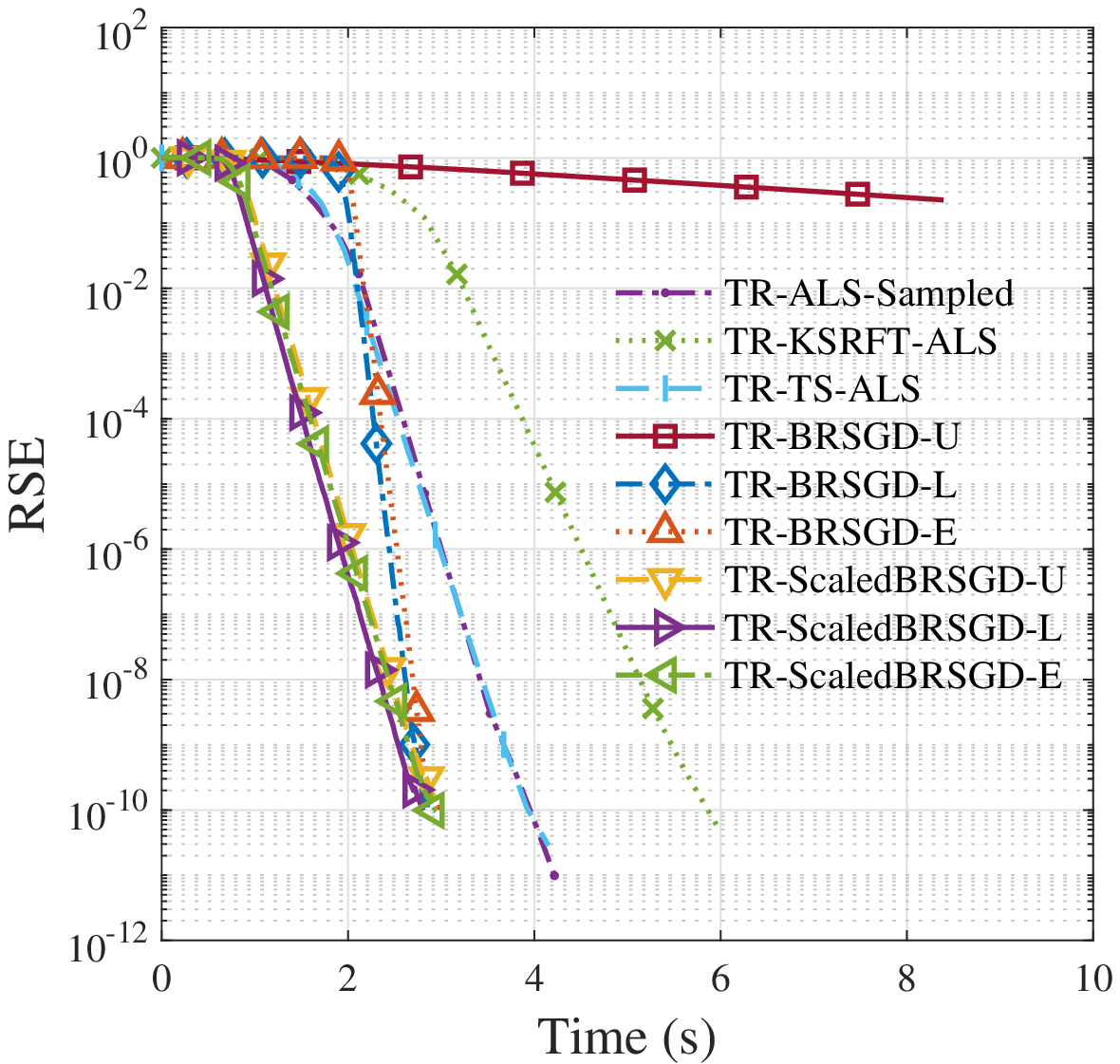}} 
	\caption{Results output by various algorithms with $tol = 1e-10$, $T = 4000$, $MT = \text{inf}$, and $R = 10$. We set the parameters as follows: TR-ALS-Sampled\&TR-KSRFT-ALS\&TR-TS-ALS: $J = 4500$; TR-GD: $\alpha = 6e-7$; TR-ScaledGD: $\alpha =  5e-1$, TR-BRSGD: $\alpha = 6e-2$, $|\fn|=200$; TR-ScaledBRSGD: $\alpha = 1e-1$, $|\fn|=200$,$|\hn|=200$.} 
	\label{fig:synA1}
\end{figure}

\begin{figure}[htbp] 
	\centering 
	\subfloat[Number of iterations v.s. RSE]{\includegraphics[scale=0.3]{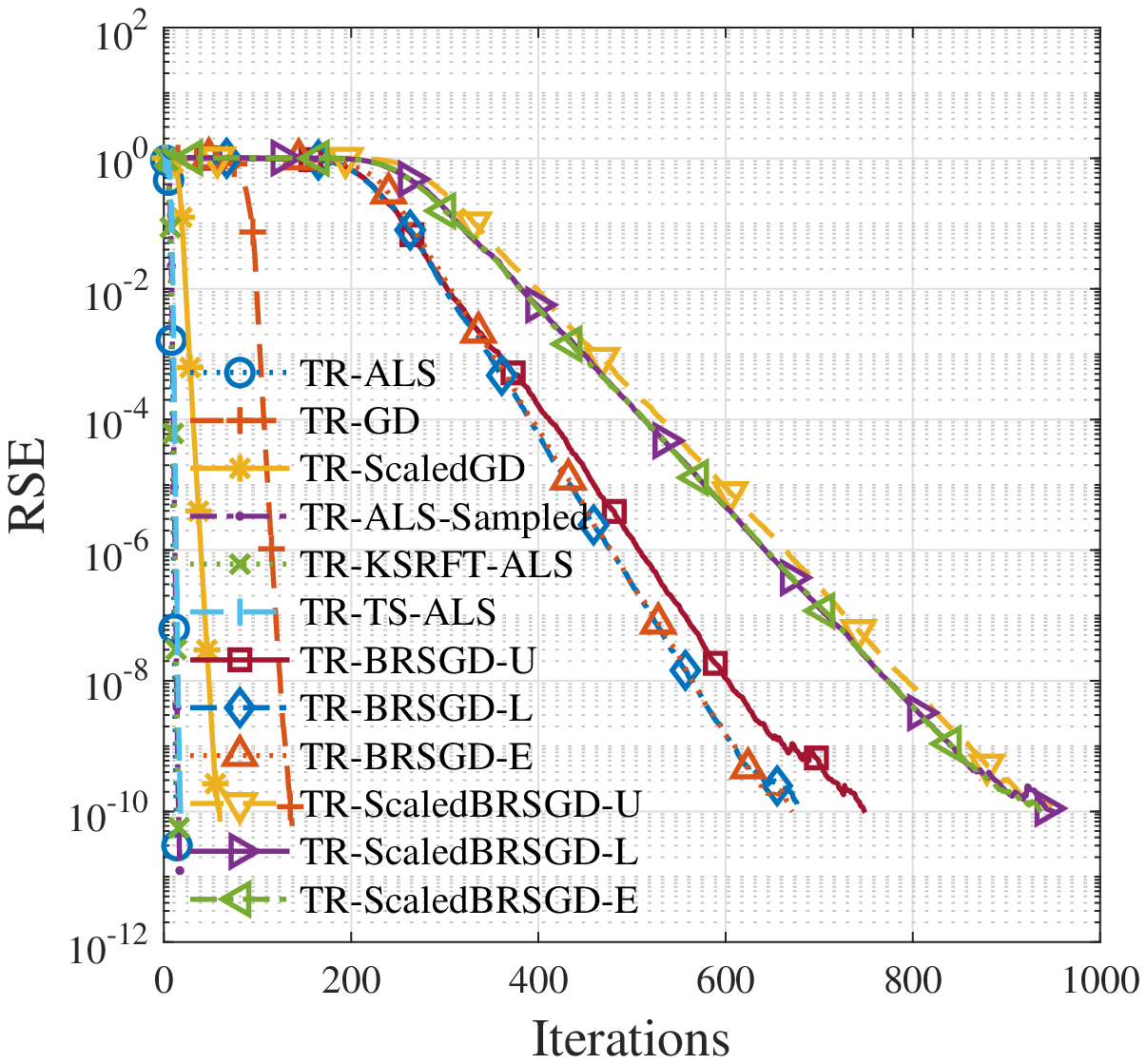}} 
	\subfloat[Time v.s. RSE]{\includegraphics[scale=0.3]{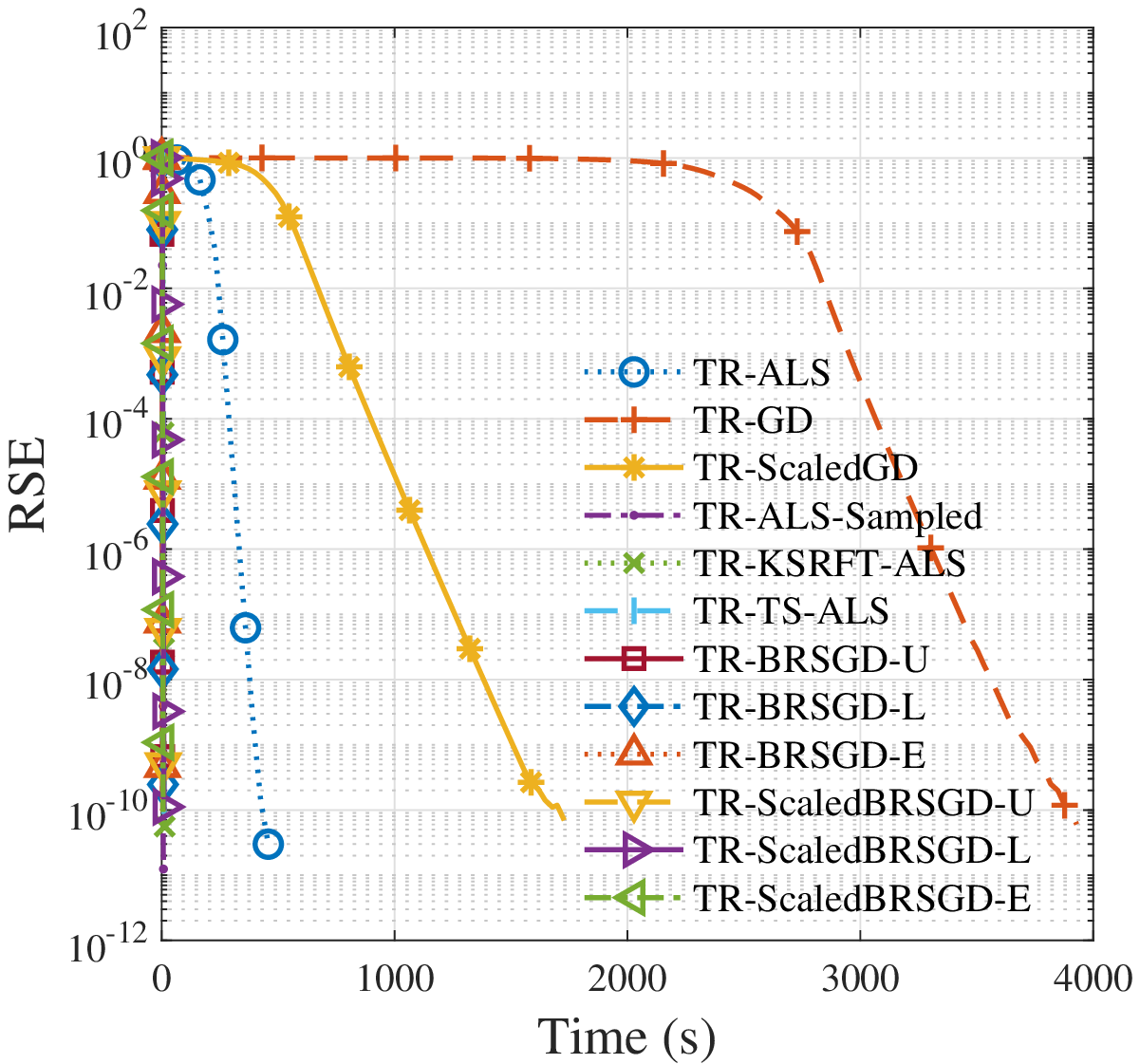}} 
	\subfloat[Zoom in Time v.s. RSE]{\includegraphics[scale=0.3]{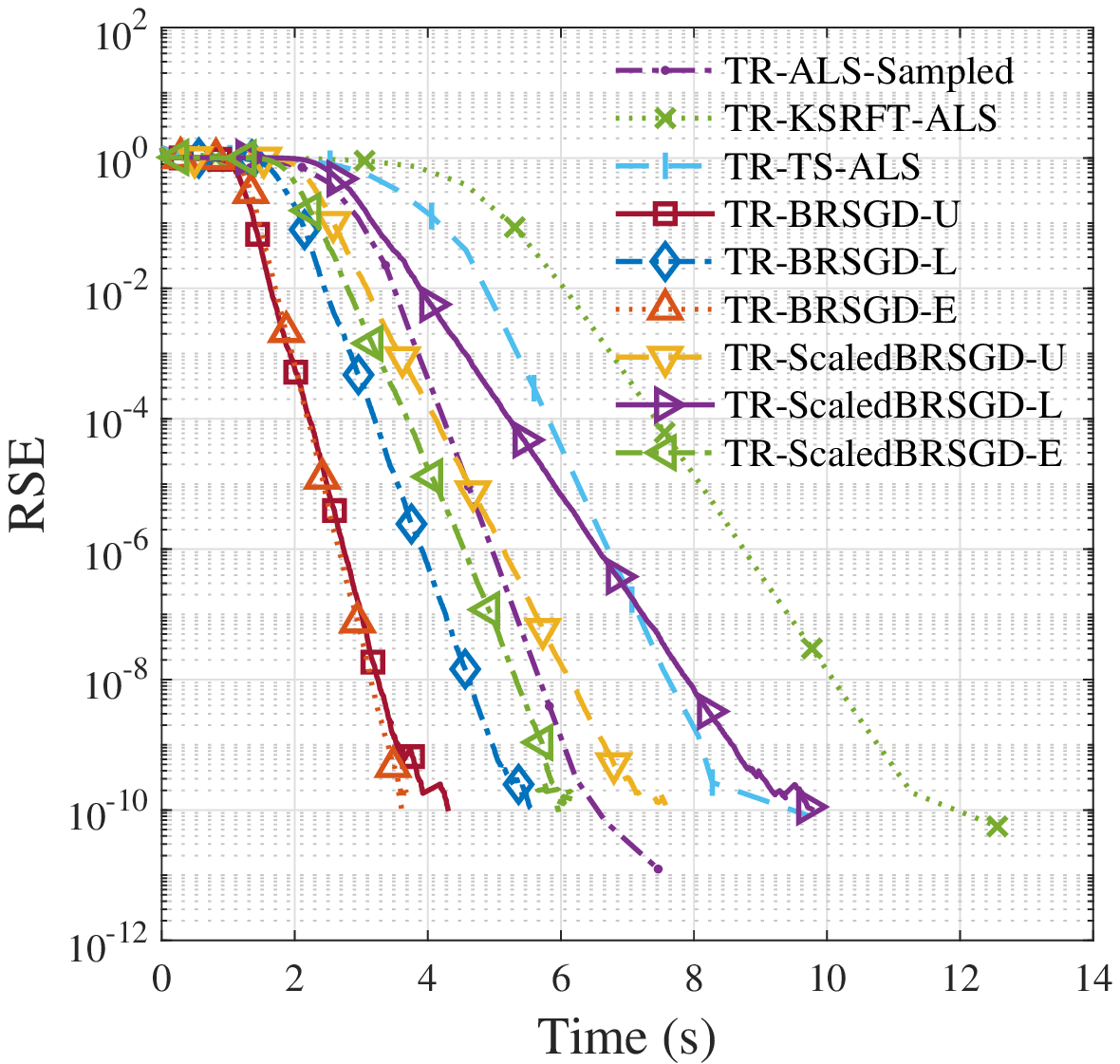}} 
	\caption{Results output by various algorithms with $tol = 1e-10$, $T = 4000$, $MT = \text{inf}$, and $R = 10$. We set the parameters as follows: TR-ALS-Sampled\&TR-KSRFT-ALS\&TR-TS-ALS: $J = 4500$; TR-GD: $\alpha = 6e-7$; TR-ScaledGD: $\alpha =  5e-1$, TR-BRSGD: $\alpha = 6e-2$, $|\fn|=200$; TR-ScaledBRSGD: $\alpha = 1e-1$, $|\fn|=200$,$|\hn|=400$.} 
	\label{fig:synA2}
\end{figure}

\paragraph{Experiment-\uppercase\expandafter{\romannumeral1}.}
The first experiment is to test our algorithms on well-conditioned data for both convergence rate and running time.
We consider $300 \times 300 \times 300$ tensors whose TR-cores of size $10 \times 300 \times 10$ are random tensors with entries drawn independently from a standard normal distribution.
The experiment results are shown in \Cref{fig:synA1,fig:synA2}. 
Some specific discussions are as follows.
\begin{itemize}
\item  For this data, all the ALS-based algorithms, i.e., TR-ALS, TR-ALS-Sampled, TR-KSRFT-ALS, and TR-TS-ALS, exhibit better convergence rates than the gradient-based algorithms, i.e., TR-GD, TR-ScsledGD, TR-BRSGD\footnote{The suffixes U, L, and E denote %the probability distributions used in the algorithms with $U$, $L$, and $E$ standing for 
the uniform, leverage-based, and  Euclidean-based probability distributions used in the algorithms, respectively.}, and TR-ScaledBRSGD, as expected. However, the total running time of the former is longer compared with the latter. This is mainly because the cost of the single-iteration in the former is expensive.

\item  Among the ALS-based algorithms, the three randomized ones, i.e., TR-ALS-Sampled, TR-KSRFT-ALS, and TR-TS-ALS, are faster than TR-ALS in both the single-iteration and total process but they spend more total time than TR-BRSGD and TR-ScaledBRSGD. This is mainly because the sketching sizes of the former are much larger than the batch sizes involved in the latter, resulting in the single-iteration of the former being slower than the one in our methods. This becomes especially obvious when the dimensionality and order of the tensor data grow.

\item  Among all the gradient-based algorithms,  as expected, TR-BRSGD and TR-ScaledBRSGD require much less overall running time than TR-GD and TR-ScaledGD.
Also, it can be seen that TR-ScaledBRSGD is faster than TR-BRSGD in \Cref{fig:synA1}, but the case is reversed in \Cref{fig:synA2}. 
This is mainly because the step size has a significant impact on the convergence speed of the stochastic algorithms. In contrast, 
TR-ScaledGD and TR-GD are less affected by the step size. In addition, it is not surprising to see that the two deterministic methods %former 
always converges faster than our methods.

\item Since this data is very uniform and well-conditioned, there is no noticeable difference in performance between different sampling probability distributions. Nevertheless, from the small difference, we can still find that the algorithms based on importance sampling converge faster. 
\end{itemize}

\paragraph{Experiment-\uppercase\expandafter{\romannumeral2}.}
Our second experiment considers the data with ill-conditioned TR-cores created by
\begin{equation*}
	\mat{G}_{n(2)} = \mat{U}_{n} \mat{S} \mat{V}^\intercal,
\end{equation*} 
where $\mat{U}_{n} \in \bb{R}^{I_n \times R_{true}^2}$ is a random matrix with orthonormal columns, $\mat{S} \in \bb{R}^{R_{true}^2 \times R_{true}^2}$ is a diagonal matrix with diagonal entries geometrically decreasing such that the condition number of $\mat{G}_{n(2)}$ is $\kappa$, and $\mat{V} \in \bb{R}^{R_{true}^2 \times R_{true}^2}$ is an orthogonal matrix shared by all TR-cores.

By varying $\kappa$, we can have some different tensors and further obtain the numerical results of this experiment, which are reported in \Cref{tab:synB}. From this table, it is seen that TR-ScaledBRSGD, especially with importance sampling, has better performance. That is, it can achieve better accuracy in the same time, which means that less time is needed for convergence. 
By contrast, TR-GD and TR-BRSGD cannot achieve satisfactory accuracy in a short time because of the slow convergence due to the ill-conditionedness of data.
Also, we can see that the sketching size needed for TR-ALS-Sampled/TR-KSRFT-ALS/TR-TS-ALS in the ill-conditioned case becomes larger, which can lead to the inefficiency of the algorithms. This may be due to the extreme sensitivity and poor sample complexity brought by the data.
Finally, it can be seen that under the same time, our two methods, i.e., TR-BRSGD and TR-ScaledBRSGD, are able to iterate more times, which further indicates that our proposed algorithms have less single-iteration time and are therefore very suitable for large-scale problems.

\begin{table}[htbp]
\centering
\caption{Results output by various algorithms with $tol = 1e-10$, $T = 1e10$, $MT = 10$, and $R = 5$. We set the parameters as follows: TR-ALS-Sampled\&TR-KSRFT-ALS\&TR-TS-ALS: $J = 8000$; TR-GD: $\alpha = 5e-6$; TR-ScaledGD: $\alpha =  8e-1$, TR-BRSGD: $\alpha = 5e-4$, $|\fn|=200$; TR-ScaledBRSGD: $\alpha = 9e-2$, $|\fn|=200$, $|\hn|=1200$.} 
\label{tab:synB}
\resizebox{1\linewidth}{!}{
\begin{tabular}{llllllllll}
    \toprule
    & \multicolumn{3}{c}{$\kappa=1e2$} & \multicolumn{3}{c}{$\kappa=1e4$} & \multicolumn{3}{c}{$\kappa=1e6$} \\
    \cmidrule(lr){2-4}
    \cmidrule(lr){5-7}
    \cmidrule(lr){8-10}
    Method 		 & {RSE} & {Iterations} & {Time} & {RSE} & {Iterations} & {Time} & {RSE} & {Iterations} & {Time}\\
    \midrule
    TR-ALS 	           & 6.91e-1 & 1.00e0 & 1.55e1   & 7.50e-1 & 1.00e0 & 1.57e1   & 7,26e-1 & 1.00e0 & 1.59e1\\
    TR-GD 	           & 4.20e2  & 2.00e0 & 1.76e1   & 4.43e2  & 2.00e0 & 1.75e1   & 4.36e2  & 2.00e0 & 1.79e1\\
    TR-ScaledGD 	   & 8.11e-1 & 2.00e0 & 1.73e1   & 8.61e-1 & 2.00e0 & 1.73e1   & 8.37e-1 & 2.00e0 & 1.80e1\\
    \midrule
    TR-ALS-Sampled     & 3.07e-4 & 1.27e1 & 1.03e1   & 1.11e-1 & 1.29e1 & 1.04e1   & 1.37e-1 & 1.27e1 & 1.04e1\\
    TR-KSRFT-ALS       & 1.64e-1 & 8.00e0 & 1.07e1   & 2.82e-1 & 8.00e0 & 1.06e1   & 2.96e-1 & 8.00e0 & 1.06e1\\
    TR-TS-ALS          & 3.13e-2 & 1.39e1 & 1.04e1   & 1.05e-1 & 1.40e1 & 1.04e1   & 8.60e-2 & 1.34e1 & 1.03e1\\
    \midrule 
    TR-BRSGD-U         & 1.63e2  & 1.92e3 & 1.00e1   & 1.79e3  & 1.90e3 & 1.00e1   & 1.91e2 & 1.81e3 & 1.00e1\\
    TR-BRSGD-E         & 1.61e2  & 1.92e3 & 1.00e1   & 1.69e3  & 1.94e3 & 1.00e1   & 1.84e2 & 1.84e3 & 1.00e1\\
    TR-BRSGD-L         & 1.75e2  & 1.81e3 & 1.00e1   & 1.80e3  & 1.83e3 & 1.00e1   & 1.98e2 & 1.75e3 & 1.00e1\\
    \midrule
    TR-ScaledBRSGD-U   & 2.84e-2 & 1.52e3 & 1.00e1   & 5.10e-2 & 1.53e3 & 1.00e1   & 6.70e-2 & 1.42e3 & 1.00e1\\
    TR-ScaledBRSGD-E   & 1.55e-7 & 1.43e3 & 1.00e1   & 5.72e-2 & 1.56e3 & 1.00e1   & 3.39e-6 & 1.47e3 & 1.00e1\\
    TR-ScaledBRSGD-L   & 2.58e-8 & 1.54e3 & 1.00e1   & 2.89e-2 & 1.45e3 & 1.00e1   & 8.79e-2 & 1.36e3 & 1.00e1\\
    \bottomrule
\end{tabular}
}
\end{table}	

%%==================
\subsection{Image experiments}
In this experiment, we test our algorithms on hyperspectral images, which are special images with two spatial coordinates and one spectral coordinate. That is, the first two orders are the image height and width, and the third one is the number of spectral bands. 
We consider three data tensors available at \url{http://www.ehu.eus/ccwintco/index.php/Hyperspectral Remote Sensing Scenes}. Their sizes are shown in \Cref{tab:hsi_result} along with our experimental results, where the parameters of different algorithms are listed in \Cref{tab:hsi_para}.

\begin{table}[htbp]
    \centering
    \caption{Results output by various algorithms for real datasets with different parameters listed in \Cref{tab:hsi_para}.}  
\label{tab:hsi_result}
\resizebox{1\linewidth}{!}{
\begin{tabular}{lllllllllllll}
    \toprule
    & \multicolumn{4}{c}{SalinasA.($83 \times 86 \times 224$)} & \multicolumn{4}{c}{IndianP.($145 \times 145 \times 220$)} & \multicolumn{4}{c}{PaviaU.($610 \times 340 \times 103$)} \\
    \cmidrule(lr){2-5}
    \cmidrule(lr){6-9}
    \cmidrule(lr){10-13}
    Method 		 & {RSE} & {PSNR} & {Iterations} & {Time}  & {RSE} & {PSNR} & {Iterations} & {Time}  & {RSE} & {PSNR} & {Iterations} & {Time}\\
    \midrule
    TR-ALS 	           & 4.33e-2 & 4.21e1 & 6.00e0 & 5.38e0   & 4.17e-2 & 3.74e1 & 2.00e0 & 1.28e1   & 2.83e-1 & 2.23e1 & 1.00e0 & 5.23e1\\
    TR-GD 	           & 1.03e-1 & 3.48e1 & 7.00e0 & 5.38e0   & 6.91e-2 & 3.28e1 & 2.00e0 & 1.09e1   & 3.38e-1 & 2.01e1 & 1.00e0 & 4.77e1\\
    TR-ScaledGD 	   & 6.54e-2 & 3.86e1 & 7.00e0 & 5.31e0   & 6.57e-2 & 3.32e1 & 2.00e0 & 1.10e1   & 3.38e-1 & 2.01e1 & 1.00e0 & 4.77e1\\
    \midrule
    TR-ALS-Sampled     & 4.52e-2 & 4.19e1 & 3.13e2 & 5.01e0   & 4.64e-2 & 3.65e1 & 1.82e2 & 1.00e1   & 1.91e-1 & 2.78e1 & 1.03e2 & 1.01e1\\
    TR-KSRFT-ALS       & 4.17e-2 & 4.23e1 & 1.63e2 & 5.02e0   & 4.53e-2 & 3.67e1 & 1.14e2 & 1.01e1   & 1.96e-1 & 2.71e1 & 4.63e1 & 1.01e1\\
    TR-TS-ALS          & 4.94e-2 & 4.10e1 & 2.45e2 & 5.01e0   & 4.67e-2 & 3.64e1 & 1.63e2 & 1.00e1   & 2.04e-1 & 2.64e1 & 9.63e1 & 1.01e1\\
    \midrule
    TR-BRSGD-U         & 1.04e-1 & 3.46e1 & 2.27e3 & 5.00e0   & 6.91e-2 & 3.28e1 & 3.15e3 & 1.00e1   & 3.42e-1 & 2.00e1 & 1.83e3 & 1.00e1\\
    TR-BRSGD-E         & 1.03e-1 & 3.48e1 & 2.24e3 & 5.00e0   & 6.91e-2 & 3.28e1 & 3.19e3 & 1.00e1   & 3.42e-1 & 1.98e1 & 1.87e3 & 1.00e1\\
    TR-BRSGD-L         & 1.04e-1 & 3.47e1 & 2.04e3 & 5.00e0   & 6.91e-2 & 3.28e1 & 1.92e3 & 1.00e1   & 3.43e-1 & 1.99e1 & 1.25e3 & 1.00e1\\
    \midrule
    TR-ScaledBRSGD-U   & 4.34e-2 & 4.20e1 & 1.48e3 & 5.00e0   & 3.82e-2 & 3.82e1 & 1.49e3 & 1.00e1   & 1.70e-1 & 2.86e1 & 1.07e3 & 1.00e1\\
    TR-ScaledBRSGD-E   & 4.30e-2 & 4.20e1 & 1.50e3 & 5.00e0   & 3.81e-2 & 3.82e1 & 1.52e3 & 1.00e1   & 1.70e-1 & 2.86e1 & 8.44e3 & 1.00e1\\
    TR-ScaledBRSGD-L   & 4.32e-2 & 4.21e1 & 1.40e3 & 5.00e0   & 3.88e-2 & 3.80e1 & 1.15e3 & 1.00e1   & 1.70e-1 & 2.85e1 & 1.04e3 & 1.00e1\\
    \bottomrule
\end{tabular}
}
\end{table}	

\begin{table}[htbp]
    \centering
    \caption{Parameters of different algorithms for real datasets.} 
    \label{tab:hsi_para}
        \resizebox{1\linewidth}{!}{
    \begin{tabular}{llll}
        \toprule
        & \multicolumn{3}{c}{$tol = 1e-10$, $T = 1e10$, $MT = 5$} \\
        \cmidrule(lr){2-4}
        Method 		 & {SalinasA.($R = 5$)} & {IndianP.($R = 10$)} & {PaviaU.($R = 10$)}\\
        \midrule
        \tabincell{l}{TR-ALS-Sampled\\ \&TR-KSRFT-ALS\\ \&TR-TS-ALS} 	           & $J=400$ & $J=800$ & $J=1000$   \\
        TR-GD           & $\alpha = 2e-13$ & $\alpha = 2e-16$ & $\alpha = 2e-16$   \\
        TR-ScaledGD     & $\alpha = 1e-1$ & $\alpha = 1e-1$ & $\alpha = 1e-1$   \\
        TR-BRSGD ($|\fn|=200$)	    & $\alpha = 2e-9$ & $\alpha = 5e-14$ & $\alpha = 8e-10$   \\
        TR-ScaledBRSGD ($|\fn|=200$,$|\hn|=1000$) 	& $\alpha = 2e-2$ & $\alpha = 4e-3$ & $\alpha = 4e-2$   \\
        \bottomrule
    \end{tabular}
     }
\end{table}	

From \Cref{tab:hsi_result}, we can observe similar results to the ones from previous experiments. That is, TR-ScaledBRSGD performs best under both the measurements RSE and PSNR in almost the same running time, which illustrates the very strong applicability of our algorithms to all types of data.
Furthermore, from the view of the variation of parameters, the batch sizes of our two algorithms do not grow with the increase of data size or TR-ranks compared to TR-ALS-Sampled/TR-KSRFT-ALS/TR-TS-ALS. This means that the cost of the single-iteration in our methods can be controlled.

%%=================%%
%%   Conclusions   %%
%%=================%%
\section{Concluding Remarks}
\label{sec:conclusion}
In this paper, we develop a doubly randomized optimization framework for large-scale TR decomposition, i.e., TR-BRSGD, and also present its scaled variant, i.e., TR-ScaledBRSGD. They make full use of the structure of TR decomposition and work well when suitable probability distributions are employed for selecting the mini-batch.
Numerical experiments are given to illustrate this conclusion. 

The presentation in this paper mainly focuses on the problem without constraints. Actually, our methods %and theory
entail considerable flexibility, and many frequent regularizers and constraints can be readily handled. 
In addition, 
the framework can also be applied to various other high-dimensional problems, such as tensor robust principal component analysis, tensor completion, tensor regression, and so on. We leave them for future research.

%%=================%%
%%  Declarations   %%
%%=================%%
\section*{Declarations}
\subsection*{Funding}
This work was supported by the National Natural Science Foundation of China (No. 11671060) and the Natural Science Foundation Project of CQ CSTC (No. cstc2019jcyj-msxmX0267).

\subsection*{Data Availability}
The data that support the findings of this study are available from the corresponding author upon reasonable request.

\subsection*{Competing Interests}
The authors declare that they have no conflict of interest.

% \begin{acknowledgements}
% If you'd like to thank anyone, place your comments here.
% \end{acknowledgements}

\bibliographystyle{spmpsci}      % mathematics and physical sciences
\bibliography{bit-bibliography}   % name your BibTeX data base

% % Non-BibTeX users please use
% \begin{thebibliography}{}
% %
% % and use \bibitem to create references. Consult the Instructions
% % for authors for reference list style.
% % Format for Journal Reference
% Author, Article title, Journal, Volume, page numbers (year)
% Author, Book title, page numbers. Publisher, place (year)
% \end{thebibliography}

\appendix\section*{Appendix}
\section{Proofs for lemmas and theorems} %Proof of \Cref{ssec:TR-BRSGD}}
\subsection{Proof of \Cref{lem:unbias_gradient}}
\label{appdxpf:unbias_gradient}
Assume that we sample the $j_f$-th gradient with the probability $p_{j_f}$, and use $|\fxi|$ gradients to estimate the full gradient at the $t$-th iteration. 
Further, for $f = 1,\cdots, |\fxi|$, define
\begin{equation*}
	\mf{G}_f = \frac{1}{|\fxi| J_{\xi^{t}} p_{j_f}} \nabla_{\mat{G}_{\xi^{t}(2)}}f_{j_f}(\tensor{Y}^{t}).
\end{equation*} 
Thus,
\begin{align*}
	\bb{E} [\mf{G}_f] &= \frac{1}{|\fxi| J_{\xi^{t}}} \sum_{j_f = 1}^{J_{\xi^{t}}} \left( p_{j_f} \frac{\nabla_{\mat{G}_{\xi^{t}(2)}}f_{j_f}(\tensor{Y}^{t})}{p_{j_f}} \right) 
	= \frac{1}{|\fxi|} \nabla_{\mat{G}_{\xi^{t}(2)}}f(\tensor{Y}^{t}).
\end{align*}
Since 
$\mf{g}_{\xi^{t}}^t = \sum_{f=1}^{|\fxi|} \mf{G}_f$, we have the desired result
\begin{equation*}
	\bb{E}_{\zeta^{t}}\left[ \mf{g}_{\xi^{t}}^t~|~\mc{B}^{t},\xi^{t} \right] = \sum_{f=1}^{|\fxi|} \bb{E} [\mf{G}_f] = \nabla_{\mat{G}_{\xi^{t}(2)}}f(\tensor{Y}^{t}).
\end{equation*}

\subsection{Proof of \Cref{thm:conv_brsgd}}
\label{appdxpf:thm:conv_brsgd}
Considering the block-wise Lipschitz continuity, for $\hat{\tensor{Y}}$, $\bar{\tensor{Y}}$ and the mode $n \in \{1,\cdots, N \}$, there exists a constant $L_{n}$ such that 
\begin{equation}
\label{eq:lip_grad}
f(\hat{\tensor{Y}}) \leq f(\bar{\tensor{Y}}) + \left< \nabla_{\mat{G}_{n(2)}} f(\bar{\tensor{Y}}), \hat{\mat{G}}_{n(2)} - \bar{\mat{G}}_{n(2)} \right> + \frac{L_n}{2} \left\| \hat{\mat{G}}_{n(2)} - \bar{\mat{G}}_{n(2)} \right\|_F^2,
\end{equation}
where $\hat{\tensor{Y}}$ and $\bar{\tensor{Y}}$ are the same except the $n$-th TR-core, i.e., $\hat{\mat{G}}_{i(2)}= \bar{\mat{G}}_{i(2)}$ for $i \neq n$. Hence, combining with \eqref{eq:update_brsgd}, we have
\begin{align}
\label{eq:lip_grad1}
    f(\tensor{Y}^{t+1})-f(\tensor{Y}^{t}) 
    &\leq \langle \nabla_{\mat{G}_{\xi^{t}(2)}} f(\tensor{Y}^{t}), \mat{G}_{\xi^{t}(2)}^{t+1}-\mat{G}_{\xi^{t}(2)}^{t} \rangle + \frac{L}{2} \left\|\mat{G}_{\xi^{t}(2)}^{t+1}-\mat{G}_{\xi^{t}(2)}^{t} \right\|_F^2\nonumber\\
    &= -\alpha^t \langle \nabla_{\mat{G}_{\xi^{t}(2)}} f(\tensor{Y}^{t}), \mf{g}_{\xi^{t}}^t \rangle + \frac{(\alpha^t)^2 L}{2} \left\| \mf{g}_{\xi^{t}}^t \right\|_F^2,
\end{align}
where $L = \max_{t=0,\cdots,\infty} L_{\xi^{t}}^t \leq \infty$.
Note that the updates $\mat{G}_{n(2)}^{t}$ are bounded for all $n$ and $t$. Then 
\begin{align}\label{eq:smooth}
\left\| \mf{g}_{\xi^{t}}^t \right\|_F^2 \leq M \textrm{ and } M < \infty. 
\end{align}
Thus, taking the expectation of the inequality \eqref{eq:lip_grad1} conditioned on the filtration $\mc{B}^{t}$ and $\xi^{t}$ shows that
\begin{align*}
	\bb{E}_{\zeta^{t}}\left[ f(\tensor{Y}^{t+1})~|~\mc{B}^{t},\xi^{t} \right] - f(\tensor{Y}^{t})
	&\leq -\alpha^t \left\| \nabla_{\mat{G}_{\xi^{t}(2)}} f(\tensor{Y}^{t}) \right\|_F^2 + \frac{(\alpha^t)^2 L}{2} \bb{E}_{\zeta^{t}}\left[ \left\| \mf{g}_{\xi^{t}}^t \right\|_F^2~|~\mc{B}^{t},\xi^{t} \right] \\
	&\leq -\alpha^t \left\| \nabla_{\mat{G}_{\xi^{t}(2)}} f(\tensor{Y}^{t}) \right\|_F^2 + \frac{(\alpha^t)^2 L M}{2}.
\end{align*}
Further, noting $\bb{E}_{\xi^{t},\zeta^{t}}\left[ \left\| \nabla_{\mat{G}_{\xi^{t}(2)}} f(\tensor{Y}^{t}) \right\|_F^2 \right] = \frac{1}{N} \left\| \nabla f(\tensor{Y}^{t}) \right\|_F^2$ and taking the total expectation of the above inequality, we get
\begin{equation*}
	\bb{E} \left[ f(\tensor{Y}^{t+1}) \right] - \bb{E} \left[ f(\tensor{Y}^{t}) \right] \leq -\frac{\alpha^t}{N} \bb{E} \left[ \left\| \nabla f(\tensor{Y}^{t}) \right\|_F^2 \right] + \frac{(\alpha^t)^2 L M}{2}.
\end{equation*}
Summing up the inequality from $t = 0$ to $t = T$ and taking $T \rightarrow \infty$, denoting $f(\tensor{Y}^*)$ as the global optimal value, and combining with $f(\tensor{Y}) \geq f(\tensor{Y}^*)$, it implies that 
\begin{equation*}
	\sum_{t=0}^{\infty} \frac{\alpha^t}{N} \bb{E} \left[ \left\| \nabla f(\tensor{Y}^{t}) \right\|_F^2 \right] \leq f(\tensor{Y}^0) - f(\tensor{Y}^*) + \sum_{t=0}^{\infty} \frac{(\alpha^t)^2 L M}{2}.
\end{equation*} 
Note that $\sum_{t=0}^{\infty} (\alpha^t)^2 < \infty$, thus, using \cite[Lemma A.5]{mairal2010OnlineLearning}, we can get the desired result
\begin{equation*}
	\liminf_{t \rightarrow \infty} \bb{E} \left[ \left\| \nabla f(\tensor{Y}^{t}) \right\|_F^2 \right] = 0.
\end{equation*}

%\section{Proof of \Cref{ssec:TR-BRSGD-probability}}
\subsection{Proof of \Cref{thm:variance}}
\label{appdxpf:variance}
Define
\begin{equation*}
\Theta_f = \frac{1}{|\fxi| J_{\xi^{t}} p_{j_f}} \nabla_{\mat{G}_{\xi^{t}(2)}} f_{j_f}(\tensor{Y}),
\end{equation*}
where $f = 1,\cdots, |\fxi|$.
Thus,
\begin{equation*}
\bb{E}[\Theta_f] = \sum_{j_f=1}^{J_{\xi^{t}}} \frac{1}{|\fxi| J_{\xi^{t}} p_{j_f}} p_{j_f} \nabla_{\mat{G}_{\xi^{t}(2)}} f_{j_f}(\tensor{Y}) = \frac{1}{|\fxi|} \nabla_{\mat{G}_{\xi^{t}(2)}}  f(\tensor{Y}),
\end{equation*}
and
\begin{align*}
\bb{V}[(\Theta_f)_{i,r}] 
&= \bb{E}[(\Theta_f)_{i,r}^2] - \bb{E}^2[(\Theta_f)_{i,r}] \\
&= \frac{1}{|\fxi|^2 J_{\xi^{t}}^2} \sum_{j_f=1}^{J_{\xi^{t}}} \left( \frac{\left[ \nabla_{\mat{G}_{\xi^{t}(2)}} f_{j_f}(\tensor{Y}) \right]_{i,r}^2}{p_{j_f}} \right) - \frac{1}{|\fxi|^2} \left[ \nabla_{\mat{G}_{\xi^{t}(2)}} f(\tensor{Y}) \right]_{i,r}^2.
\end{align*}
Since $\bb{V} \left[ \left( \mf{g}_{\xi^{t}}^{t} \right)_{i,r} \right] = \sum_{f=1}^{|\fxi|} \bb{V}[(\Theta_f)_{i,r}]$, we have
\begin{equation*}
\bb{V} \left[ \left( \mf{g}_{\xi^{t}}^{t} \right)_{i,r} \right] = \frac{1}{|\fxi| J_{\xi^{t}}^2} \sum_{j_f=1}^{J_{\xi^{t}}} \left( \frac{\left[ \nabla_{\mat{G}_{\xi^{t}(2)}} f_{j_f}(\tensor{Y}) \right]_{i,r}^2}{p_{j_f}} \right) - \frac{1}{|\fxi|} \left[ \nabla_{\mat{G}_{\xi^{t}(2)}} f(\tensor{Y}) \right]_{i,r}^2.
\end{equation*}
On the other hand, 
\begin{align*}
&\bb{E}_{\zeta^{t}}\left[ \left\| \mf{g}_{\xi^{t}}^{t} - \nabla_{\mat{G}_{\xi^{t}(2)}} f(\tensor{Y}^{t}) \right\|_F^2 ~|~ \mc{B}^{t},\xi^{t} \right] \\
&\qquad = \sum_{i=1}^{I_{\xi^{t}}} \sum_{r=1}^{R_{\xi^{t}} R_{\xi^{t}+1}} \bb{E} \left[ \left( \mf{g}_{\xi^{t}}^{t} - \nabla_{\mat{G}_{\xi^{t}(2)}} f(\tensor{Y}) \right)_{i,r}^2 \right]  = \sum_{i=1}^{I_{\xi^{t}}} \sum_{r=1}^{R_{\xi^{t}} R_{\xi^{t}+1}} \bb{V} \left[ \left( \mf{g}_{\xi^{t}}^{t} \right)_{i,r} \right].
\end{align*}
Therefore,
\begin{align*}
&\bb{E}_{\zeta^{t}}\left[ \left\| \mf{g}_{\xi^{t}}^{t} - \nabla_{\mat{G}_{\xi^{t}(2)}} f(\tensor{Y}^{t}) \right\|_F^2 ~|~ \mc{B}^{t},\xi^{t} \right] \\
&\quad = \frac{1}{|\fxi| J_{\xi^{t}}^2} \sum_{j_f=1}^{J_{\xi^{t}}} \frac{1}{p_{j_f}} \sum_{i=1}^{I_{\xi^{t}}} \sum_{r=1}^{R_{\xi^{t}} R_{\xi^{t}+1}} \left[ \nabla_{\mat{G}_{\xi^{t}(2)}} f_{j_f}(\tensor{Y}) \right]_{i,r}^2 - \frac{1}{|\fxi|} \left\| \nabla_{\mat{G}_{\xi^{t}(2)}} f(\tensor{Y}) \right\|_F^2 \\
&\quad = \frac{1}{|\fxi|} \sum_{j_f=1}^{J_{\xi^{t}}} \frac{1}{p_{j_f}} \left( \sum_{i=1}^{I_{\xi^{t}}} (\mat{R}_{[\xi^{t}]}^{t})_{i,j_f}^2 \right) \left( \sum_{r=1}^{R_{\xi^{t}} R_{\xi^{t}+1}} ((\mat{G}^{\ne \xi^{t}}_{[2]})^t)_{j_f,r}^2 \right) \\
&\qquad - \frac{1}{|\fxi|} \left\| \nabla_{\mat{G}_{\xi^{t}(2)}} f(\tensor{Y}) \right\|_F^2 \\
&\quad = \frac{1}{|\fxi|} \sum_{j_f=1}^{J_{\xi^{t}}} \frac{1}{p_{j_f}} \left\| \mat{R}_{[\xi^{t}]}^{t}(:,j_f) \right\|_2^2 \left\| (\mat{G}^{\ne \xi^{t}}_{[2]})^t (j_f,:) \right\|_2^2 - \frac{1}{|\fxi|} \left\| \nabla_{\mat{G}_{\xi^{t}(2)}} f(\tensor{Y}) \right\|_F^2,
\end{align*} 	
where the second equality follows from the result $$\nabla_{\mat{G}_{\xi^{t}}} f_{j_f}(\tensor{Y}) = J_{\xi^{t}} \mat{R}_{[\xi^{t}]}^{t}(:,j_f) (\mat{G}^{\ne \xi^{t}}_{[2]})^t (:,j_f).$$ 

To show that  the probability distribution $\vect{q}^{\ne \xi^{t}}$ given in \eqref{eq:optimal_prob} minimizes $\bb{E}_{\zeta^{t}}\left[ \left\| \mf{g}_{\xi^{t}}^{t} - \nabla_{\mat{G}_{\xi^{t}(2)}} f(\tensor{Y}^{t}) \right\|_F^2 ~|~ \mc{B}^{t},\xi^{t} \right]$,
we first define a function as
\begin{equation*}
f(p_1, \cdots p_{J_{\xi^{t}}}) = \sum_{j_f=1}^{J_{\xi^{t}}} \frac{1}{p_{j_f}} \left\| \mat{R}_{[\xi^{t}]}^{t}(:,j_f) \right\|_2^2 \left\| (\mat{G}^{\ne \xi^{t}}_{[2]})^t (j_f,:) \right\|_2^2,
\end{equation*}
which characterizes the dependence of the variance 
\begin{equation*}
\bb{E}_{\zeta^{t}}\left[ \left\| \mf{g}_{\xi^{t}}^{t} - \nabla_{\mat{G}_{\xi^{t}(2)}} f(\tensor{Y}^{t}) \right\|_F^2 ~|~\mc{B}^{t},\xi^{t} \right]
\end{equation*}
on the sampling probability distribution $\vect{p}=[p_1, \cdots, p_{J_{\xi^{t}}}]^\intercal$. %To minimize $f$ subject to $\sum_{i=1}^{J_{\xi^{t}}} p_i = 1$, 
Further, we introduce the Lagrange multiplier $\lambda$ and define the following function
\begin{equation*}
g(p_1, \cdots p_{J_{\xi^{t}}}) = f(p_1, \cdots p_{J_{\xi^{t}}}) + \lambda \left( \sum_{i=1}^{J_{\xi^{t}}} p_i - 1 \right).
\end{equation*}
Since
\begin{equation*}
0 = \frac{\partial g}{\partial p_i} = \frac{-1}{p_i^2} \left\| \mat{R}_{[\xi^{t}]}^{t}(:,i) \right\|_2^2 \left\| (\mat{G}^{\ne \xi^{t}}_{[2]})^t(i,:) \right\|_2^2 + \lambda,
\end{equation*}
we have
\begin{equation*}
p_i = \frac{\left\| \mat{R}_{[\xi^{t}]}^{t}(:,i) \right\|_2^2 \left\| (\mat{G}^{\ne \xi^{t}}_{[2]})^t(i,:) \right\|_2^2}{\sqrt{\lambda}} 
= \frac{\left\| \mat{R}_{[\xi^{t}]}^{t}(:,i) \right\|_2^2 \left\| (\mat{G}^{\ne \xi^{t}}_{[2]})^t(i,:) \right\|_2^2}{\sum_{i'=1}^{J_n} \left\| \mat{R}_{[\xi^{t}]}^{t}(:,i') \right\|_2^2 \left\| (\mat{G}^{\ne \xi^{t}}_{[2]})^t(i',:) \right\|_2^2},
\end{equation*}
where the second equality is from the fact that $\sum_{i=1}^{J_{\xi^{t}}} p_i = 1$. Note that, for the above probabilities, $\frac{\partial^2 g}{\partial p_i^2} > 0$. 
Hence, the probability distribution in \eqref{eq:optimal_prob} indeed minimizes the variance. Meanwhile,
 substituting $\vect{q}^{\ne \xi^{t}}$ in \eqref{eq:optimal_prob} into \eqref{eq:var} gives \eqref{eq:optimal_var}.
So, all the desired results hold.
	
\subsection{Proof of \Cref{lem:euc_est}}
\label{appdxpf:euc_est}
Clearly, all $\vect{q}^{\ne n}(i)$ are nonnegative. Moreover,
\begin{equation*}
	\sum_{\substack{i_{n+1}, \cdots, i_N \\ i_{1}, \cdots, i_{n-1}}} \vect{q}^{\ne n}(\overline{i_{n+1} \cdots i_{N} i_{1} \cdots i_{n-1}}) = \prod_{\substack{k=1\\k \ne n}}^{N} \sum_{i_k = 1}^{I_k} \vect{p}_k(i_k) = \prod_{\substack{k=1\\k \ne n}}^{N} 1 = 1,
\end{equation*}
due to the fact that each $\vect{p}_k$ is a probability distribution. So, $\vect{q}^{\ne n}$ is also a probability distribution. 
Furthermore, for each $n \in [N]$, we have
\begin{align*}
    \vect{q}^{\ne n}(\overline{i_{n+1} \cdots i_{N} i_{1} \cdots i_{n-1}}) &= \prod_{\substack{k=1\\k \ne n}}^{N} \vect{p}_k(i_k) 
	= \frac{\prod_{k \ne n} \| \tensor{G}_{k}(:,i_k,:) \|_F^2}{\prod_{k \ne n} \| \tensor{G}_{k} \|_F^2} \\
    &\geq \frac{\| \tensor{G}_{n+1}(:,i_{n+1},:) \cdots \tensor{G}_{N}(:,i_{N},:) \tensor{G}_{1}(:,i_{1},:) \cdots \tensor{G}_{n-1}(:,i_{n-1},:) \|_F^2}{\prod_{k \ne n} \| \tensor{G}_{k} \|_F^2} \\
    &= \frac{\| \tensor{G}^{\ne n}(:,i,:) \|_F^2}{\prod_{k \ne n} \| \tensor{G}_{k} \|_F^2} 
	= \beta_n \| \tensor{G}^{\ne n}(:,i,:) \|_F^2
	\geq \beta_n \vect{p}^{\ne n}(i)
\end{align*}
as desired.

%\section{Proof of \Cref{ssec:TR-BRSGD-Scaled}}
\subsection{Proof of \Cref{lem:unbias_hessian}}
\label{appdxpf:unbias_hessian}
Assume that we sample the $j_h$-th gradient with probability $p_{j_h}$, and use $|\hxi|$ gradients to estimate the full gradient at the $t$-th iteration. 
Further, for $h = 1,\cdots, |\hxi|$, define
\begin{equation*}
	\mf{H}_h = \frac{1}{|\fxi| J_{\xi^{t}} p_{j_h}} \nabla^2_{{\rm vec}(\mat{G}_{\xi^{t}(2)}),{\rm vec}(\mat{G}_{\xi^{t}(2)})}f_{j_h}(\tensor{y}_{t}).
\end{equation*} 
Thus,
\begin{align*}
	\bb{E} [\mf{H}_h] &= \frac{1}{|\hxi| J_{\xi^{t}}} \sum_{j_h = 1}^{J_{\xi^{t}}} \left( p_{j_h} \frac{\nabla^2_{{\rm vec}(\mat{G}_{\xi^{t}(2)}),{\rm vec}(\mat{G}_{\xi^{t}(2)})}f_{j_h}(\tensor{Y}_{t})}{p_{j_h}} \right) \\
	&= \frac{1}{|\hxi|} \nabla^2_{{\rm vec}(\mat{G}_{\xi^{t}(2)}),{\rm vec}(\mat{G}_{\xi^{t}(2)})}f(\tensor{Y}_{t}),
\end{align*}
which together with $\tilde{\mf{h}}_{\xi^{t}}^t = \sum_{h=1}^{|\hxi|} \mf{H}_h$ 
implies the desired result:
\begin{equation*}
	\bb{E}_{\zeta^{t}}\left[ \tilde{\mf{h}}_{\xi^{t}}^t~|~\mc{B}^{t},\xi^{t} \right] = \sum_{h=1}^{|\hn|} \bb{E} [\mf{H}_h] = \nabla^2_{{\rm vec}(\mat{G}_{\xi^{t}(2)}),{\rm vec}(\mat{G}_{\xi^{t}(2)})}f(\tensor{Y}^{t}),
\end{equation*}
and hence
\begin{equation}
\label{eq:h_unbias}
    \bb{E}_{\zeta^{t}}\left[ \mf{h}_{\xi^{t}}^t~|~\mc{B}^{t},\xi^{t} \right] = (\mat{G}^{\ne \xi^{t}}_{[2]})^\intercal \mat{G}^{\ne \xi^{t}}_{[2]}.
    \end{equation}

\subsection{Proof of \Cref{lem:unbias_sd}}
\label{appdxpf:unbias_sd}
From our assumption, we have $\left\| \nabla^2_{{\rm vec}(\mat{G}_{\xi^{t}(2)}),{\rm vec}(\mat{G}_{\xi^{t}(2)})}f(\tensor{Y}^t) \right\| \leq 1$, and hence $\left\| \mf{h}_{\xi^{t}}^t \right\| \leq 1$, where $\| \cdot \|$ denotes the spectral norm of a  matrix. Note that for a positive definite matrix $\mat{A}$ satisfying $\|\mat{A}\| \leq 1$, 
\begin{equation*}
%\label{eq:inv_taylor}
	\mat{A}^{-1} = \sum_{i=0}^{\infty}(\mat{I}-\mat{A})^{i}.
\end{equation*}
Then, 
\begin{equation}\label{eq:h_taylor}
	(\mf{h}_{\xi^{t}}^t)^{-1}= \sum_{i=0}^{\infty} (\mat{I} - \mf{h}_{\xi^{t}}^t)^i \approx \mat{I}+\mat{I}-\mf{h}_{\xi^{t}}^t = 2\mat{I}-\mf{h}_{\xi^{t}}^t.
\end{equation}
On the other hand, using \eqref{eq:h_unbias}, we obtain
\begin{equation*}
	\bb{E}_{\zeta^{t}}\left[ 2\mat{I}-\mf{h}_{\xi^{t}}^t ~|~\mc{B}^{t},\xi^{t} \right] 
	= 2\mat{I}-(\mat{G}^{\ne \xi^{t}}_{[2]})^\intercal \mat{G}^{\ne \xi^{t}}_{[2]}
	\approx \left( (\mat{G}^{\ne \xi^{t}}_{[2]})^\intercal \mat{G}^{\ne \xi^{t}}_{[2]} \right)^{-1}.
\end{equation*}
Thus, by the fact that the sampled sets $\fn$ and $\hn$ are independent in our algorithm and the property of expectation, we can have the desired result.

\subsection{Proof of \Cref{thm:conv_scaled}}
\label{appdxpf:conv_scaled}
%Recall that $\xi^{t}$ and $\zeta^{t}$ are the random mode and fibers chosen at the iteration $t$, respectively. Thus, under the assumption that the updates $\mat{G}_{\xi^{t}(2)}^{t}$ are bounded for all $n$ and $t$, 
Combining \eqref{eq:lip_grad}, \eqref{eq:update_brsgd}, and \eqref{eq:scaledsg_direction}, we observe:
\begin{align*}
	f(\tensor{Y}^{t+1})-f(\tensor{Y}^{t}) 
	&\leq \langle \nabla_{\mat{G}_{\xi^{t}(2)}}f(\tensor{Y}^t), \mat{G}_{\xi^{t}(2)}^{t+1} - \mat{G}_{\xi^{t}(2)}^{t} \rangle+\frac{L}{2} \left\| \mat{G}_{\xi^{t}(2)}^{t+1} - \mat{G}_{\xi^{t}(2)}^{t} \right\|_F^2 \\
	&= -\alpha^t \langle \nabla_{\mat{G}_{\xi^{t}(2)}}f(\tensor{Y}^t), \mf{g}_{\xi^{t}}^t (\mf{h}_{\xi^{t}}^t)^{-1} \rangle + \frac{(\alpha^t)^2 L}{2} \left\| \mf{g}_{\xi^{t}}^t (\mf{h}_{\xi^{t}}^t)^{-1} \right\|_F^2,
\end{align*}
where $L = \max_{t=0,\cdots,\infty} \leq L_{\xi^{t}}^t < \infty$.
% Further, define the weighted norm as $$\left\|\nabla_{\mat{G}_{\xi^{t}(2)}}f(\tensor{Y}^t) \right\|_{\mat{G}^{\ne \xi^{t}}_{[2]}}^* = \left\| \nabla_{\mat{G}_{\xi^{t}(2)}}f(\tensor{Y}^t) \left( (\mat{G}^{\ne \xi^{t}}_{[2]})^\intercal \mat{G}^{\ne \xi^{t}}_{[2]} \right)^{-1/2} \right\|_F.$$ 
Thus, taking expectation conditioned on the filtration $\mc{B}^{t}$ and the chosen mode index $\xi^{t}$, and considering \Cref{lem:unbias_sd} implies
\begin{align}
\label{eq:scaled_exp}
    \bb{E}_{\zeta^{t}} \left[ f(\tensor{Y}^{t+1}) ~|~\mc{B}^{t},\xi^{t} \right] -f(\tensor{Y}^{t}) &\lesssim -\alpha^t \left\| \nabla_{\mat{G}_{\xi^{t}(2)}}f(\tensor{Y}^t) \left( (\mat{G}^{\ne \xi^{t}}_{[2]})^\intercal \mat{G}^{\ne \xi^{t}}_{[2]} \right)^{-1/2} \right\|_F^2\nonumber \\
    &\quad+ \frac{(\alpha^t)^2 L}{2} \bb{E}_{\zeta^{t}} \left[ \left\| \mf{g}_{\xi^{t}}^t (\mf{h}_{\xi^{t}}^t)^{-1} \right\|_F^2 ~|~\mc{B}^{t},\xi^{t} \right] \nonumber\\
    &\leq \text{------}%~~~~\leq -\alpha^t \left\| \nabla_{\mat{G}_{\xi^{t}(2)}}f(\tensor{Y}^t) \left( (\mat{G}^{\ne \xi^{t}}_{[2]})^\intercal \mat{G}^{\ne \xi^{t}}_{[2]} \right)^{-1/2} \right\|_F^2 \nonumber\\
    %&\qquad
    + \frac{(\alpha^t)^2 L}{2} \bb{E}_{\zeta^{t}} \left[ \left\| \mf{g}_{\xi^{t}}^t \right\|_F^2 \left\| (\mf{h}_{\xi^{t}}^t)^{-1} \right\|^2 ~|~\mc{B}^{t},\xi^{t} \right]\nonumber\\
    &\leq \text{------} + \frac{(\alpha^t)^2 LM}{2} \bb{E}_{\zeta^{t}} \left[ \left\| (\mf{h}_{\xi^{t}}^t)^{-1} \right\|^2 ~|~\mc{B}^{t},\xi^{t} \right] \nonumber\\
    &= \text{------} + \frac{(\alpha^t)^2 L M}{2} \bb{E}_{\zeta^{t}} \left[ \left\| \sum_{i=0}^{\infty} (\mat{I} - \mf{h}_{\xi^{t}}^t)^i \right\|^2 ~|~\mc{B}^{t},\xi^{t} \right]\nonumber\\
    &\leq -\alpha^t \left\| \nabla_{\mat{G}_{\xi^{t}(2)}}f(\tensor{Y}^t) \left( (\mat{G}^{\ne \xi^{t}}_{[2]})^\intercal \mat{G}^{\ne \xi^{t}}_{[2]} \right)^{-1/2} \right\|_F^2 + \frac{(\alpha^t)^2 L M K}{2}, 
\end{align}
where %$\lesssim$ means less than or approximately equal, 
the third inequality is due to \eqref{eq:smooth}, % \mf{g}_{\xi^{t}}^t \right\|_F^2 \leq M \textrm{ with } M < \infty$,
the equality is based on \eqref{eq:h_taylor}, and the last inequality follows from % $\left\| \mf{h}_{\xi^{t}}^t \right\| \leq 1$ 
the result 
\begin{align*}
    \left\| \sum_{i=0}^{\infty} (\mat{I} - \mf{h}_{\xi^{t}}^t)^i \right\|^2 
    &\leq \left[ \sum_{i=0}^{\infty} \left(1 - \sqrt{\lambda_{\min}\left((\mf{h}_{\xi^{t}}^t)^\intercal \mf{h}_{\xi^{t}}^t\right)} \right)^i \right]^2 \\
    &= \frac{1}{\lambda_{\min}\left((\mf{h}_{\xi^{t}}^t)^\intercal \mf{h}_{\xi^{t}}^t\right)}
    \leq K 
    \leq \infty.
\end{align*}
Further, note that 
%{\color{red} 
\begin{equation*}
    \bb{E}_{\xi^{t},\zeta^{t}}\left[ \left\| \nabla_{\mat{G}_{\xi^{t}(2)}}f(\tensor{Y}^t) \left( (\mat{G}^{\ne \xi^{t}}_{[2]})^\intercal \mat{G}^{\ne \xi^{t}}_{[2]} \right)^{-1/2} \right\|_F^2 \right] = \frac{1}{N} \left\| (\mat{H}^{t})^{-1/2} \nabla f(\tensor{Y}^{t}) \right\|_F^2.
\end{equation*}
%}
Thus, taking the total expectation on \eqref{eq:scaled_exp}, we can get
\begin{equation*}
	\bb{E} \left[ f(\tensor{Y}^{t+1}) \right] - \bb{E} \left[ f(\tensor{Y}^{t}) \right] \leq - \frac{\alpha^t}{N} \bb{E} \left[ \left\| (\mat{H}^{t})^{-1/2} \nabla f(\tensor{Y}^{t}) \right\|_F^2 \right] + \frac{(\alpha^t)^2 L M K}{2}.
\end{equation*}
Summing up the above inequality from $t = 0$ to $t = T$ and taking $T \rightarrow \infty$, denoting $f(\tensor{Y}^*)$ as the global optimal value, and combining with $f(\tensor{Y}) \geq f(\tensor{Y}^*)$, we have 
\begin{equation*}
	\sum_{t=0}^{\infty} \frac{\alpha^t}{N} \bb{E} \left[ \left\| (\mat{H}^{t})^{-1/2} \nabla f(\tensor{Y}^{t}) \right\|_F^2 \right] \leq f(\tensor{Y}^0) - f(\tensor{Y}^*) + \sum_{t=0}^{\infty} \frac{(\alpha^t)^2 L M K}{2},
\end{equation*} 
which together with %Note 
the fact that 
$\sum_{t=0}^{\infty} (\alpha^t)^2 < \infty$ and % thus, using 
\cite[Lemma A.5]{mairal2010OnlineLearning} implies the desired result
%{\color{red}
\begin{equation*}
	\liminf_{t \rightarrow \infty} \bb{E} \left[ \left\| (\mat{H}^{t})^{-1/2} \nabla f(\tensor{Y}^{t}) \right\|_F^2 \right] = 0.
\end{equation*}
%}

\end{sloppypar}
\end{document}